\documentclass[11pt,a4paper,twoside]{amsart}

\usepackage{amsmath,amssymb,amsfonts,amsthm,mathrsfs}
\usepackage{times,hyperref,color}
\usepackage[toc,page]{appendix}

\newcommand{\zero}[1]{\underline{#1}_{\circ}}
\let\mathcal\mathscr

\newtheorem{The}{Theorem}[section]

\newtheorem{Theorem}{Theorem}[section]

\newtheorem{Proposition}[The]{Proposition}
\newtheorem{Conjecture}[The]{Conjecture}
\newtheorem{Lemma}[The]{Lemma}
\newtheorem{Corollary}[The]{Corollary}

\theoremstyle{definition}
\newtheorem{Definition}[The]{Definition}
\newtheorem{Convention}[The]{Convention}
\newtheorem{Remark}[The]{Remark}


\begin{document}


\title{
Biholomorphic equivalence to
\\
 totally nondegenerate model CR manifolds
\\
 and Beloshapka's maximum conjecture
}

\author{Masoud Sabzevari}
\address{Department of Mathematics,
University of Shahrekord, 88186-34141 Shahrekord, IRAN and School of
Mathematics, Institute for Research in Fundamental Sciences (IPM), P.
O. Box: 19395-5746, Tehran, IRAN}
\email{sabzevari@math.iut.ac.ir}

\date{\number\year-\number\month-\number\day}

\maketitle

\begin{abstract}
 Applying \'{E}lie Cartan's classical method, we show that the biholomorphic equivalence problem to a totally nondegenerate Beloshapka's model of CR dimension one and codimension $k> 1$, whence of real dimension $2+k$, is reducible to some absolute parallelism, namely to an $\{e\}$-structure on a certain prolonged manifold of real dimension either $3+k$ or $4+k$. The proof relies on the weight analysis of the  structure equations associated with the mentioned problem of equivalence. Thanks to the achieved results, we prove Beloshapka's maximum conjecture about the rigidity of his CR models of certain lengths equal or greater than three: {\sf "}\,CR automorphism Lie groups of these models do not contain any nonlinear map, preserving the origin\,{\sf "}. Here, we mainly deal with CR models of the fixed CR dimension one though the results seem generalizable by means of certain analogous proofs.
\end{abstract}

\pagestyle{headings} \markright{Biholomorphic equivalence to totally nondegenerate CR models}

\section{Introduction}
\label{section-introduction}

The notion of totally nondegenerate CR manifolds has a close connection with the theory of free Lie algebras. In order to explain this connection in CR dimension one, let $h_1$ and $h_2$ be two linearly independent elements of a certain vector space over the field $\mathbb C$. By definition (\cite{Merker-Porten-2006, Merker-Sabzevari-CEJM, Reutenauer}), the {\sl rank two complex free Lie algebra} $\mathcal{ F}$ is the smallest non-commutative and non-associative $\mathbb C$-algebra having $h_1$ and $h_2$ as its elements, with bilinear multiplication $(h,h')\mapsto[h, h']
\in
\mathcal{F},
$ satisfying the skew-symmetry and Jacobi-like identity:
\[
\aligned
0
&=
[h,\,h']
+
[h',\,h],
\\
0
&=
\big[h,\,[h',\,h'']\big]
+
\big[h'',\,[h,\,h']\big]
+
\big[h',\,[h'',\,h]\big],
\endaligned
\]
for arbitrary elements $h,\,h'\,h''\,\in\,\mathcal{F}$. Such an algebra $\mathcal{ F}$ is unique up to isomorphism. Importantly, {\it no linear
relation exists between iterated multiplications}, {\it i.e.} between
iterated Lie brackets, {\it except those generated
only by skew-symmetry and Jacobi identity}: this is the {\sl freeness} of the algebra.

Then, the elements of $\mathcal F$, designated as {\sl words}, can be rewritten as iterated Lie brackets between the letters $h_1$ and $h_2$. For instance:
\[
\big[[h_1,\,h_2],\,h_1\big],
\ \ \ \ \ \ \ \ \
\big[ h_1,\,\big[h_1,\,[h_2,\,h_1]\big]\big],
\ \ \ \ \ \ \ \ \
\big[ [h_2,h_1],\,\big[ h_1,\,[h_2,\,h_1]\big]\big].
\]
 We define the {\sl length} of each word to be the number of $h_1, h_2$ elements in it. Define $\mathcal{F}_1$ to be the $\mathbb C$-vector space generated by $h_1, h_2$ and for $\ell \geqslant 2$, let $\mathcal{F}_\ell$ be the $\mathbb C$-vector space generated by all words of the lengths $\leqslant \ell$. Then clearly $\mathcal{F} = \bigcup_{ \ell \geqslant 1}\, \mathcal{F}_\ell$ and we have the following {\it filtration}:
\begin{equation}
\label{filtration}
\mathcal F_1\subset\mathcal F_2\subset\mathcal F_3\subset\cdots
\end{equation}
on $\mathcal F$. Also let us denote by ${\sf n}_\ell$, the dimension of the complex vector space $\mathcal F_\ell$  which can be computed by means of the recursive relation, introduced in \cite[Theorem 2.6]{Merker-Porten-2006}.
By an induction based on the Jacobi identity,
it follows that each length $\ell$ word can be expressed as a linear combination of some specific words of the form:
\[
\big[
h_{i_1},\,\big[h_{i_2},\,\big[\dots\big[h_{i_{\ell-1}},\,h_{i_{\ell}}
\big]\dots\big]\big],
\]
which are called {\sl simple words} of the length $\ell$. Hence for each $\ell>1$ we have $\mathcal{F}_\ell:=\mathcal F_{\ell-1}+[\mathcal F_1, \mathcal F_{\ell-1}]$. The collection of all simple words generates
$\mathcal{F}$ as a vector space over $\mathbb C$, though actually it is not a basis\footnote{One basis for the free Lie algebra $\mathcal F$ is the so-called {\sl Hall-Witt basis} ({\it see} \cite[Definition 2.5 and Theorem 2.6]{Merker-Porten-2006}).}.

Let $\mathfrak F_\ell:=\mathcal F_\ell\setminus\mathcal F_{\ell-1}$ be the $\mathbb C$-vector space generated by all (simple) words of the length $\ell$. Then, $\mathfrak F_\ell$ is of dimension ${\sf m}_\ell:={\sf n}_\ell-{\sf n}_{\ell-1}$ and $[\frak F_{\ell_1}, \frak F_{\ell_2}]\subset\frak F_{\ell_1+\ell_2}$ for each $\ell_1, \ell_2\geqslant 1$. Consequently, our infinite dimensional free algebra $\mathcal F$ {\it is graded} of the form:
\begin{equation}
\label{gradation}
\frak F_1\oplus\frak F_2\oplus\frak F_3\oplus\ldots.
\end{equation}

Now, let us turn our attention to the subject of totally nondegenerate CR manifolds of CR dimension one. Let $M\subset\mathbb C^{1+k}$ be a real analytic generic submanifod of codimension $k$, and hence of real dimension $2+k$. As is known (\cite{BER, Boggess-1991, Merker-Porten-2006, 5-cubic}), the holomorphic subbundle $T^{1,0}M$ of its complexified bundle $\mathbb C\otimes TM$ can be generated by a single holomorphic vector field $\mathcal L$. Let $D_1:=T^{1,0}M+ T^{0,1}M$, where $T^{0,1}M:=\overline {T^{1,0}M}$ and let us also define successively $D_j=D_{j-1}+[D_1,D_{j-1}]$ for $j>1$. As is customary in the Lie-Cartan theory, we assume {\sl strong uniformity}, that is: for each $j\geqslant 1$, the dimension of $D_j$ is fully constant on the points of $M$ if it is thought of as being local. So, all $D_j$s are subbundles of
$\mathbb C\otimes TM$. It is also natural to assume that $M$ is {\em minimal} (\cite{BER,
Merker-Porten-2006, 5-cubic}), in the sense that:
\[
D_i
=
\mathbb C\otimes TM
\ \ \ \ \
\text{\rm for all}\ \
i\geqslant i^*\ \
\text{\rm large enough}.
\]
Lastly, as a first step in the study of such differential structures, it is also natural to assume that the ranks of the subbundles $D_1, D_2, D_3, \dots$
increase as much as possible.

\begin{Definition}
\label{totally-nondegenerate}
An arbitrary (local) real analytic CR generic submanifold $M\subset\mathbb C^{1+k}$ of CR dimension one and codimension $k$ is called {\sl totally nondegenerate}\,\,---\,\,or {\sl completely nondegenerate} or {\sl maximally minimal} ({\it cf.} \cite{5-cubic})\,\,---\,\,whenever $\mathbb C\otimes TM$ can be generated by means of the minimum possible number of iterated Lie brackets between the generators $\mathcal L$ and $\overline{\mathcal L}$ of $D_1$, increasing maximally through a filtration:
\[
D_1\varsubsetneq D_2\varsubsetneq\ldots\varsubsetneq D_\rho=\mathbb C\otimes TM.
\]
In this case, the length $\rho$ of this filtration is also called the {\sl length} of $M$.
\end{Definition}

Set $h_1:=\mathcal L$ and $h_2:=\overline{\mathcal L}$. Rephrasing this in the language of free Lie algebras, a real analytic CR generic submanifold $M\subset\mathbb C^{1+k}$ of CR dimension one is totally nondegenerate of the length $\rho$ whenever for each $\ell=1,\ldots,\rho-1$, the vector space $D_\ell$ and its basis can be identified by its corresponding $\mathcal F_\ell$. More precisely, whenever the rank of $D_\ell$ is maximum, equal to the dimension ${\sf n}_\ell$ of $\mathcal F_\ell$ and $D_\ell$ behaves precisely as an ${\sf n}_{\ell}$-dimensional $\mathbb C$-vector space generated by the maximum possible number $\sum_{l\leqslant\ell}\,{\sf m}_{l}$ of (simple) iterated Lie brackets between $\mathcal L$ and $\overline{\mathcal L}$ of the lengths $l\leq\ell$. As is the case with the free Lie algebra $\mathcal F$, {\it no linear relation exists between the iterated brackets of $\mathcal L$ and $\overline{\mathcal L}$ in the lengths $\leq\rho-1$, except those generated by skew-symmetry and Jacobi identity}. However, the case of the last bundle $D_\rho=\mathbb C\otimes TM$ is in part a different matter. Careful inspection of the above definition shows that the length $\rho$ of a $k$-codimensional submanifold $M$ is indeed the smallest integer $\ell$ satisfying:
\[
{\rm rank}\,\mathbb C \otimes TM\leq {\rm dim}\,\mathcal F_\ell
\]
or equivalently $2+k\leqslant {\sf n}_\ell$. It is of course possible in certain codimensions $k$ that $\rho$ satisfies $2+k\lneqq {\sf n}_\rho$. Such a constraint on the rank of the complexified bundle $\mathbb C\otimes TM$ in comparison with the maximal freedom ${\sf n}_\rho$ may cause an encounter with a length $\rho$ simple iterated bracket $\mathcal T_\rho\in D_\rho \setminus D_{\rho-1}$ which is not independent of the other length $\rho$ simple brackets {\it even} modulo skew-symmetry and Jacobi identity. In this case, $\mathcal T_\rho$ is expressible as $\mathcal T_\rho:=\sum_i\,{\sf f}_i\,\mathcal T_{\ell,i}$ for some basis simple Lie brackets $\mathcal T_{\ell,i}$ of the lengths $\ell\leqslant\rho$ and some certain functions ${\sf f}_i$ defined on $M$. This occurs for example in the case of totally nondegenerate CR manifolds $M\subset\mathbb C^3$ of CR dimension $1$ and codimension $2$, explained in \cite[p. 87]{Samuel}. Notice that a similar phenomenon can also occur for iterated Lie brackets of lengths bigger than $\rho$ ({\it see e.g.} subsections 4.4 and 4.5 of \cite{5-cubic}).

\subsection{Existence of totally nondegenerate CR manifolds in arbitrary codimensions}

For a CR generic submanifold $M\subset\mathbb C^{1+k}$ of codimension $k$, consider the {\sl complex tangent space} $T^cM:={\rm Re}\,(T^{1,0}M)$. If $\mathcal L:={\sf X}+i\,\sf Y$ is the single generator of $T^{1,0}M$, then clearly ${\sf X, Y}$ generate $T^cM$. As is known \cite{BER, Boggess-1991, Merker-Porten-2006}:
\[
\mathbb C\otimes T^cM=T^{1,0}M\oplus T^{0,1}M.
\]
Set ${\sf D}_1:=T^cM$ and define, as above, the subbundles ${\sf D}_j:={\sf D}_1+[{\sf D}_1,{\sf D}_{j-1}], j\geqslant 1$ of the real bundle $TM$. Thanks to the equality $D_1=\mathbb C\otimes {\sf D}_1$ and $[\mathbb C\otimes {\sf D}_{1}, \mathbb C\otimes {\sf D}_{j-1}]=\mathbb C\otimes [{\sf D}_{1}, {\sf D}_{j-1}]$, one verifies  by means of a simple induction that $D_j=\mathbb C\otimes {\sf D}_j$ for each $j\geqslant 1$. Therefore, it is possible to restate the definition of total nondegeneracy in terms of the real distributions ${\sf D}_j$ as follows\footnote{For technical reasons, we prefer to keep this definition in terms of the complexified distributions $D_j$ as Definition \ref{totally-nondegenerate}.}: $M$ is a length $\rho$ totally nondegenerate submanifold of $\mathbb C^{1+k}$ if and only if its associated tangent bundle $TM$ can be generated by means of the minimum possible number of iterated Lie brackets between the generators ${\sf X}:={\rm Re}\,\mathcal L$ and ${\sf Y}:={\rm Im}\, \mathcal L$ of ${\sf D}_1=T^cM$, increasing maximally through the filtration:
\[
{\sf D}_1\varsubsetneq {\sf D}_2\varsubsetneq\ldots\varsubsetneq {\sf D}_\rho= TM.
\]

Roughly speaking, $M$ is a length $\rho$ totally nondegenerate CR manifold whenever after setting $h_1:=\sf X$ and $h_2:=\sf Y$, the behavior of the above filtration and also the Lie brackets between $\sf X, Y$ can be identified by the (simple) words of the lengths $\leq\rho$, belonging to the rank two ({\it real}) free Lie algebra.

Now, let us consider the existence of totally nondegenerate CR manifolds, a question which may arise naturally at this time. For a fixed positive integer $k$\,\,---\,\,which will take over the role of the codimension for the sought CR manifolds\,\,---\,\,let $\rho$ be the smallest length $\ell$ such that ${\sf n}_\ell:={\rm dim}\,\mathcal F_\ell\geqslant 2+k$. Then by \cite[Theorem 2.7]{Merker-Porten-2006}, one can find a rank two real subdistribution ${\sf D}_1:=\langle {\sf X}, {\sf Y} \rangle$ of $T\,\mathbb R^{2+k}$, defined on a neighborhood $\Omega\subset\mathbb R^{2+k}$ of the origin such that:
\begin{itemize}
\item[(i)] ${\rm dim}\,{\sf D}_\ell(0)={\sf n}_\ell$ for each $\ell<\rho$ and
\item[(ii)] ${\rm dim}\,{\sf D}_\rho(0)=2+k$ or equivalently ${\sf D}_\rho(0)=T_0\,\mathbb R^{2+k}$,
\end{itemize}
where, as above, ${\sf D}_j:={\sf D}_{j-1}+[{\sf D}_1, {\sf D}_{j-1}]$ and where ${\sf D}_\ell(0)\subset T_0\mathbb R^{2+k}$ is the image space of ${\sf D}_j$ at the origin. By definition (\cite[p. 74]{BER}), the second property (ii) indicates that the origin $0\in\mathbb R^{2+k}$ is a {\it finite type} point of the distribution ${\sf D}_1$ of the full type $2+k$. Since ${\sf D}_\rho$ has the maximum possible dimension at $0$, one finds a certain open subset $M$ of $\Omega$, including the origin, such that ${\sf D}_\rho (p)$ is again of the maximum possible dimension for each $p\in M$, {\it i.e.} ${\sf D}_\rho(p)=T_p\mathbb R^{2+k}$. Then $M$, as an open subset of $\mathbb R^{2+k}$, is a real submanifold of dimension $2+k$ and we claim that it is actually a CR generic submanifold of $\mathbb C^{1+k}\equiv\mathbb R^{2+2k}$ of codimension $k$. In fact, we define the complex structure map $J:{\sf D}_1\rightarrow {\sf D}_1$ on the subbundle ${\sf D}_1:={\sf D}_1\mid_M$ of $ TM$ by $J({\sf X})={\sf Y}$ and $J({\sf Y})=-{\sf X}$. Accordingly, the complexified bundle $\mathbb C\otimes {\sf D}_1$ can be decomposed as ({\it cf.} \cite[p. 1573]{Isaev-Zaitsev-2013}):
\[
\mathbb C\otimes {\sf D}_1:=D^{1,0}\oplus D^{0,1}
\]
where $D^{1,0}$ and $D^{0,1}$ are generated by the single vector fields $\mathcal L:={\sf X}+i\,{\sf Y}$ and $\overline{\mathcal L}:={\sf X}-i\,{\sf Y}$. The rank one complex subbundle $D^{1,0}\subset\mathbb C\otimes TM$ is involutive since clearly $[D^{1,0},D^{1,0}]\subset D^{1,0}$. Then by definition (\cite{Isaev-Zaitsev-2013}), ${\sf D}_1$ is a {\it CR structure} and the real submanifold $M$ {\it is a CR manifold} of CR dimension $1=\frac{1}{2}\,{\rm rank} \, {\sf D}_1$ and codimension $k={\rm dim}_{\mathbb R}\,M-2\,{\rm CRdim}\,M$. As before, we can denote ${\sf D}_1$, $D^{1,0}$ and $D^{0,1}$ by $T^cM$, $T^{1,0}M$ and $T^{0,1}M$, respectively. Moreover, as $M$ is an open subset of $\mathbb R^{2+k}$, then it is real analytic. Hence according to \cite[Proposition 3.3]{5-cubic} (expanded version) we can assume that $M$, regarded locally, is a {\it generic} submanifold  of $\mathbb C^{1+k}$. Finally, the two properties (i) and (ii) guarantee that this $M$ is also totally nondegenerate.

By definition, on the other hand, for every arbitrary real analytic totally nondegenerate submanifold $M\subset\mathbb C^{1+k}$ of codimension $k$, passing through the origin and with ${\sf D}_1=T^cM$, the above two items (i) and (ii) are satisfied. This indicates that the so-called {\sl H\"{o}rmander numbers} of $M$ ({\it see} \cite{BER} for definition) are $2,3,4,\ldots,\rho$ with the {\it maximum} possible {\sl multiplicities} ${\sf m}_2, {\sf m}_3, \ldots, {\sf m}_{\rho-1}, {\sf m}_\rho'$, respectively, where, ${\sf m}_j:={\sf n}_j-{\sf n}_{j-1}$ for $j=1,\ldots,\rho-1$\,\,---\,\, as was in the case of the rank two free Lie algebra\,\,---\,\,and ${\sf m}_\rho':=k-\sum_{j=1}^{\rho-1}\,{\sf m}_j\leqslant {\sf m}_\rho$. Then, by summing up the results and applying Theorems 4.3.2 and 4.5.1 of \cite{BER}, we can state that;

\begin{Theorem}
\label{theorem-existence-totally-momdeg}
\begin{itemize}
\item[(i)] In each codimension $k$, totally nondegenerate real analytic CR generic submanifolds $M\subset\mathbb C^{1+k}$, passing through the origin, exist. The length $\rho$ of such manifolds is determined as the smallest integer $\ell$ such that ${\sf n}_\ell\geqslant 2+k$.
 \item[(ii)]
 In agreement with the above notations, consider the canonical coordinates $(z, {\bf w}^2, \ldots, {\bf w}^{\rho-1},  {\bf w}^\rho)$ of  $\mathbb C^{1+k}$, with $z\in\mathbb C$, with ${\bf w}^j\in\mathbb C^{{\sf m}_j}$ for $j=2,\ldots,\rho-1$ and with ${\bf w}^{\rho}\in\mathbb C^{{\sf m}_\rho'}$. Assign the weight $1$ to $z$ and the weight $j$ to each component of the vector ${\bf w}^j$ and to its real and imaginary parts, as well, for $j=1,\ldots,\rho$. Then, the already mentioned submanifolds $M\subset\mathbb C^{1+k}$ can be represented locally near the origin as the graph of some $k$ real analytic functions:
\begin{equation}
\label{totally-nondegenerate-equations}
{\rm Im}\,{\bf w}^j:={\bf \Phi}_j(z,\overline z, {\rm Re}\,{\bf w}^2, \ldots, {\rm Re}\,{\bf w}^{j-1})+{\rm O}(j) \ \ \ \ \ \  {\scriptstyle (j\,=\,2\,,\,\ldots\,,\,\rho)},
\end{equation}
where ${\bf \Phi}_j$ is a weighted homogeneous vector-valued polynomial of the weight $j$ and where ${\rm O}(j)$ is some certain (possibly vanishing) sum of monomials of the weights $\geqslant j+1$. Moreover, denoting by ${\bf \Xi}_j:={\bf\Phi}_j+{\rm O}(j)$, the right hand sides of the above equations, we have:
\[
{\bf\Xi}_j(0,\overline z, {\rm Re}\,{\bf w})\equiv {\bf\Xi}_j(z,0, {\rm Re}\,{\bf w})\equiv 0.
\]
\end{itemize}
\end{Theorem}

\medskip

In 2004, Valerii Beloshapka established in \cite{Beloshapka2004} his universal model surfaces associated with totally nondegenerate CR manifolds and designed an effective method for their construction. It was actually along the celebrated approach initiated first by Henri Poincar\'{e} \cite{Poincare} in 1907 to study real submanifolds in the complex space $\mathbb C^2$ by means of the associated model surface, namely the {\it Heisenberg sphere} (\cite{Merker-Sabzevari-CEJM}). Several years later in 1974, Chern and Moser in their seminal work \cite{Chern-Moser} notably developed this approach by associating appropriate models to nondegenerate real {\it hypersurfaces} in complex spaces. In this framework, many questions about automorphism groups, classification, invariants and others can be reduced to similar problems about the associated models.

To the best of the author's knowledge, Beloshapka's work is the most general model-construction in the class of totally nondegenerate CR manifolds of arbitrary dimensions. Roughly speaking and after appropriate weight assignments ({\it see} $\S$\ref{section-Beloshapkas-models} for more details), a Beloshapka's model for the class of all CR manifolds in CR dimension one (as is our specific case), codimension $k$ and accordingly determined length $\rho$ is represented as the graph of some weighted homogeneous polynomial functions of the form:
 \begin{equation*}
{\rm Im}\,{\bf w}^j:={\bf \Phi}_j(z,\overline z, {\rm Re}\,{\bf w}^2, \ldots, {\rm Re}\,{\bf w}^{j-1}) \ \ \ \ \ \  {\scriptstyle (j\,=\,2\,,\,\ldots\,,\,\rho)},
\end{equation*}
obtained actually by removing the non-homogeneous parts ${\rm O}(j)$ from the general defining equations \thetag{\ref{totally-nondegenerate-equations}} of manifolds belonging to this class. Beloshapka's models are all homogeneous, of finite type and enjoy several other {\it nice} properties (\cite[Theorem 14]{Beloshapka2004}) that exhibit their significance. Two totally nondegenerate CR manifolds are holomorphically equivalent whenever their associated models are as well. Moreover, they are most symmetric nondegenerate surfaces in the sense that the dimension of the group of automorphisms ({\it see} below for definition) associated with a totally nondegenerate  manifold does not exceed that of its model.

\begin{Convention}
 Let us stress that throughout this paper, we mainly deal with Beloshapka's totally nondegenerate CR generic models in CR dimension one which, for the sake of brevity, are also termed {\sl "CR models"} or {\sl "models"}. We fix the notation $M_k$ for such CR models in codimension $k$.
\end{Convention}

For a length $\rho$ CR model $M_k\subset\mathbb C^{1+k}$ in coordinates $(z,w_1,\ldots,w_k)$, a holomorphic vector field:
\[
{\sf X}:=Z(z,w)\frac{\partial}{\partial z}+\sum_{l=1}^k\,W^l(z,w)\frac{\partial}{\partial w_l}
\]
is called an {\sl infinitesimal CR automorphism} whenever its real part is tangent to $M_k$, that is $({\sf X}+\overline{\sf X})|_{M_k}\equiv 0$. The collection of all infinitesimal CR automorphisms associated with $M_k$ form a Lie algebra, denoted by $\frak{aut}_{CR}(M_k)$. It is the CR symmetry Lie algebra of $M_k$ in the terminology of Sophus Lie's symmetry theory \cite{Merker-Engel-Lie} and is of finite dimension, of polynomial type and graded\,\,---\,\,in the sense of Tanaka\,\,---\,\,of the form (\cite{Beloshapka2004, Filomat}):
\begin{equation}
\label{aut}
\frak{aut}_{CR}(M_k):=\underbrace{\frak g_{-\rho}\oplus\cdots\frak g_{-1}}_{\frak g_-}\oplus\frak g_{0}\oplus\underbrace{\frak g_{1}\oplus\cdots\oplus\frak g_{\varrho}}_{\frak g_+}, \ \ \ \ \ \ \varrho,\,\rho\,\in\,\mathbb N,
\end{equation}
with $[\frak g_i,\frak g_j]\subset\frak g_{i+j}$.
Viewing the real analytic CR generic model $M_k$ in a purely intrinsic way, one may consider the local Lie group ${\sf Aut}_{CR}(M_k)$, associated with $\frak{aut}_{CR}(M_k)$, comprising automorphisms of the CR structure, namely of local $\mathcal C^\infty$ diffeomorphisms $h : M_k\rightarrow  M_k$ satisfying:
\[
h_\ast(T^cM_k)=T^cM_k.
\]
In other words, $h$ belongs to ${\sf Aut}_{CR}(M_k)$ if and only if it is a (local) biholomorphism of $M_k$ (\cite{Merker-Pocchiola-Sabzevari-5-CR-II}). Corresponding to \thetag{\ref{aut}}, one may write:
\begin{equation}
\label{Aut}
{\sf Aut}_{CR}(M_k):={\sf G}_-\oplus {\sf G_0}\oplus {\sf G}_+.
\end{equation}
Beloshapka in \cite{Beloshapka2004} showed that the Lie group ${\sf G}_-$ associated with the above subalgebra $\frak g_-$ of $\frak{aut}_{CR}(M_k)$ is $(2+k)$-dimensional, acts on $M_k$ freely and can naturally be identified with $M_k$, itself. Also, ${\sf G}_0$ associated with the subalgebra $\frak g_0$ comprises all {\it linear} automorphisms of $M_k$ in the isotropy subgroup ${\sf Aut}_0(M_k)$ of ${\sf Aut}_{CR}(M_k)$ at $0\in M_k$  while ${\sf G}_+$, associated with $\frak g_+$, comprises as well all the {\it nonlinear} ones.

Determining such Lie algebras of infinitesimal CR automorphisms is a question which lies pivotally at the heart of the problem of classifying local analytic CR manifolds up to biholomorphisms ({\it see e.g.} \cite{Beloshapka-Kossovskiy-2011} and the references therein). In fact, the groundbreaking works of Sophus Lie and his followers show that the most fundamental question in concern here is to draw up lists of possible such Lie algebras which would classify all possible manifolds according to their CR symmetries. Moreover, having in hand these algebras also enables one to treat the problem of constructing (canonical) Cartan geometries on certain classes of CR manifolds (\cite{Merker-Sabzevari-CEJM, ERA}) or to construct the so-called moduli spaces of model real submanifolds (\cite{SCM-Moduli}).

From a computational point of view, although computing the nonpositive part $\frak g_-\oplus\frak g_0$ of $\frak{aut}_{CR}(M_k)$ is partly convenient\,\,---\,\,in particular by means of the algorithm designed in \cite{Filomat}\,\,---\,\,unfortunately for $\frak g_+$ one needs highly complicated computations which rely on constructing and solving certain arising systems of partial differential equations (\cite{Mamai2009, 5-cubic, SCM, Merker-Sabzevari-M3-C2}). Nevertheless, after several years of experience in computing these algebras in various dimensions,  Beloshapka in \cite{Beloshapka} conjectured that\footnote{Although Beloshapka introduced his conjecture in 2012 but he and his students had been aware of it since several years before ({\it see e.g.} \cite{Beloshapka2007, Gammel-Kossovskiy-2006, Mamai2009}).};

\begin{Conjecture}
\label{conjecture}
{\bf [Beloshapka's Maximum Conjecture]} Each of Beloshapka's totally nondegenerate CR models $M$ of the length $\rho\geq 3$ has {\sl rigidity}; that is: in its associated graded Lie algebra $\frak{aut}_{CR}(M)$, the subalgebra $\frak g_+$ is trivial or equivalently $\varrho=0$ ({\it cf.}  \thetag{\ref{aut}}).
\end{Conjecture}

Holding this conjecture true may bring about having several other facts about CR models or their associated totally nondegenerate CR manifolds ({\it see e.g.} \cite{Beloshapka2007}). At present, there are only a few considerable results that verify this conjecture in some specific cases. For instance, Gammel and Kossovskiy  (\cite{Gammel-Kossovskiy-2006}) confirmed it in the specific length $\rho=3$. Also, Mamai in \cite{Mamai2009} proved this conjecture for the models of CR dimension one and codimensions $k\leq 13$. Some more relevant (partial) results in this setting are also as follows:
 \begin{itemize}
\item[$\bullet$] If $\rho=2$, then $\varrho \leqslant 2$ (\cite[p. 32]{Beloshapka1996}).
\item[$\bullet$] If $\rho=4$, then $\varrho \leqslant 1$ (\cite[Corollary 7]{Beloshapka2001}).
\item[$\bullet$] If $\rho=5$, then $\varrho \leqslant k$, where $k$ is the CR codimension
of $M$ (\cite[Proposition 2.2]{Shananina2004}).
\end{itemize}
In almost all of these works the results are achieved by means of {\it directly} computing  the associated desired Lie algebras. But the difficulty of this method which lies in the incredible {\it differential-algebraic} complexity involved ({\it cf.} \cite{SCM}), {\it may convince one of the necessity of attacking this conjecture through other ways}.

On the other hand, recently in \cite{5-cubic} and in particular in $\S 5$ of this paper, we attempted to study, by means of Cartan's classical approach, the biholomorphic equivalence problem to the $5$-dimensional length $3$ cubic model $M_3\subset\mathbb C^4$\,\,---\,\,denoted there by $M^5_{\sf c}$\,\,---\,\,represented as the graph of three defining equations:
\[
\aligned
{\rm Im}\,w_1 = z\overline{z}, \ \ \ \ \ \
{\rm Im}\,w_2  =z\overline{z}(z+\overline{z}), \ \ \ \ \ \
{\rm Im}\,w_3  = z\overline{z}(z-\overline{z}).
\endaligned
\]
As can be observed ({\it cf.} \cite[Theorem 5.1]{5-cubic}), the associated $7$-dimensional CR automorphism algebra:
\[
\frak{aut}_{CR}(M_3):=\frak g_{-3}\oplus\frak g_{-2}\oplus\frak g_{-1}\oplus\frak g_0,
 \]
 computed in $\S 3$ of this paper, is surprisingly isomorphic to that defined by the final {\sl constant type} structure equations of the mentioned equivalence problem to $M_3$. This observation was our original motivation to look upon Cartan's classical approach as an appropriate way to study Beloshapka's maximum conjecture. Examining this idea on some other CR models like those studied in \cite{Merker-Sabzevari-CEJM, Samuel, SCM-Moduli,  Merker-Sabzevari-M3-C2} also convinced us more about the effectiveness of this approach to suit our purpose. Indeed, the systematic approach developed in recent years by Jo\"el Merker, Samuel Pocchiola and the present author provides a unified way toward treating the wide variety of biholomorphic equivalence problems between CR manifolds.

 Cartan's classical method for solving equivalence problems includes three major parts: {\it absorption, normalization} and {\it prolongation}. In the CR context, usually all steps require advanced computations, the size of which increases considerably as soon as the dimension of CR manifolds increases, even by one unit. In particular, among the absorbtion-normalization steps, one encounters some arising polynomial systems, the solutions of which determine the value of some group parameters associated with the problem.
As is quite predictable due to the arbitrariness of dimension, one of our main obstacles to proving Beloshapka's maximum conjecture by applying Cartan's classical approach is actually solving these arising polynomial systems in this general manner, namely the outcome of normalizing the group parameters. In order to bypass and manipulate this critical complexity, our main weapon in this paper is in fact some helpful results achieved by a careful {\it weight analysis} of the equivalence problems, under study. Such analysis enables us to provide a much more convenient weighted homogeneous subsystem of the already mentioned system which is {\it deceptively hidden} inside the original one and opens our way of finding the desired general outcome of the normalization process.

This paper is organized as follows. In the next preliminary section, Section \ref{section-Beloshapkas-models}, we present a brief description of constructing defining equations of Beloshapka's CR models in CR dimension one.

 Next in Section \ref{section-lifted-coframes}, we endeavor to find certain expression of the so-called structure equations associated with the biholomorphic equivalence problem between an arbitrary length $\rho$ CR model $M_k$ and any arbitrary totally nondegenerate CR manifold ${\bf M}_k$ of the same length and codimension ({\it cf.} Theorem \ref{theorem-existence-totally-momdeg}). For this purpose, first we construct, in an almost explicit manner, an initial frame:
\[
\mathbb L:=\big\{\mathcal L_{1,1}, \mathcal L_{1,2}, \mathcal L_{2,3}, \ldots, \mathcal L_{\rho,2+k}\big\}
 \]
 on $\mathbb C\otimes TM_k$, where $\mathcal L_{1,1}$ and $\mathcal L_{1,2}=\overline{\mathcal L_{1,1}}$ are the single generators of $T^{1,0}M_k$ and $T^{0,1}M_k$, respectively and where each $\mathcal L_{\ell,i}$ is a length $\ell$ iterated Lie bracket between them {\it constructed as a simple word}. We also consider a so-called lifted frame ${\bf L}:=\{{\bf L}_{1,1}, {\bf L}_{1,2}, \ldots, {\bf L}_{\rho,2+k}\}$ on $\mathbb C\otimes T{\bf M}_k$, constructed as simple words written by the single generators ${\bf L}_{1,1}$ and ${\bf L}_{1,2}:=\overline{{\bf L}_{1,1}}$. Let $\Sigma:=\{\sigma_{1,1}, \sigma_{1,2}, \ldots, \sigma_{\rho,2+k}\}$ and $\Gamma:=\{\Gamma_{1,1}, \Gamma_{1,2}, \ldots, \Gamma_{\rho,2+k}\}$ be  two coframes on $\mathbb C\otimes T^\ast M_k$  and $\mathbb C\otimes T^\ast {\bf M}_k$, dual to the frames $\mathbb L$ and $\bf L$, respectively. We realize that for a general biholomorphic equivalence map $h:M_k\rightarrow {\bf M}_k$, the associated matrix of the induced complexified linear pull-back $h^\ast:\mathbb C\otimes T^\ast{\bf M}_k\rightarrow \mathbb C\otimes T^\ast M_k$,
  expressed in terms of the coframes $\Gamma$ and $\Sigma$ is an invertible $(2+k)\times (2+k)$ lower triangular matrix:
\[
\footnotesize\aligned
{\bf g}:=\left(
  \begin{array}{ccccccc}
    a_1^{p}\overline a_1^{q} & \cdot\cdot &  0 & 0 & 0 & 0 \\
    \vdots & \vdots &  \vdots & \vdots & \vdots & \vdots \\
    a_\bullet &  \ldots &  a_3 & a_1\overline a_1 & 0 & 0 \\
    a_\bullet  & \ldots &  a_5 & -\overline a_2 & \overline a_1 & 0 \\
    a_\bullet & \ldots &  a_4 & a_2 & 0 & a_1 \\
  \end{array}
\right) \ \ \ \ \ \ \ \ \ {\rm with} \ \ \ \ \ a_1\neq 0,
\endaligned
\]
for some certain complex-valued functions $a_1, a_2, a_3, \ldots$, in terms of the coordinates $(z,\overline z, w, \overline w)$ of $\mathbb C^{1+k}$. Only some powers of $a_1$ and $\overline a_1$ are visible at the diagonal of $\bf g$. As is standard in the terminology of Cartan's theory, we call this matrix the {\it ambiguity matrix} (or {\it G-structure}) of the mentioned equivalence problem and its nonzero entries $a_\bullet$ the {\it group parameters}. The collection of all such matrices forms a Lie group $G$ which is called the {\it structure Lie group} of the equivalence problem.

The main focus of Section \ref{section-weight-analysis} is on a weight analysis on the structure equations, constructed in the preceding section by applying necessary differentiations and computations on the already obtained equality $\Gamma={\bf g}\cdot\Sigma$. In particular, after appropriate weight assignment to the appearing group parameters and also after inspecting carefully the inverse of the ambiguity matrix $\bf g$, we discover that all the torsion coefficients appearing through the structure equations are of the same weight zero (Proposition \ref{prop-torsion-weight-zero}).

Next in Section \ref{section-picking-system}, we consider the outcome of the absorption and normalization steps on the constructed structure equations. It is in this section that we extract a subtle weighted homogeneous subsystem of the polynomial system, arising among the absorbtion and normalization steps. Solving this subsystem by means of some computational techniques from {\sl weighted algebraic geometry} (\cite{Alg-Geometry}), we conclude that except $a_1$, all the other group parameters (or functions) $a_2, a_3, \ldots$ must be vanishing;

 \begin{Proposition}
 (cf. Proposition \ref{prop-aj=0})
All the appearing group parameters $a_2, a_3, \ldots$ vanish identically after sufficient steps of absorption and normalization.
  \end{Proposition}

This key result converts our structure equations into a delicate constant type form ({\it cf.} Proposition \ref{prop-struc-after-vanish}):
\begin{equation}
\label{int-struc-eq}
\aligned
d\,\Gamma_{\ell,i}&:=(p_i\,\alpha+q_i\,\overline\alpha)\wedge\Gamma_{\ell,i}
+\sum_{l+m=\ell\atop{j,n}}\,{\sf c}^{i}_{j,n}\,\Gamma_{l,j}\wedge\Gamma_{m,n} \ \ \ \ {\scriptstyle (\ell\,=\,1\,,\,\ldots\,,\,\rho\,,\, \ \ i\,=\,1\,,\,\ldots\,,\,2+k).}
\endaligned
\end{equation}
  Here, $p_i, q_i, {\sf c}^i_{j,n}$ are some constant integers and $\alpha:=\frac{da_1}{a_1}$ is the only remaining Maurer-Cartan 1-form.
 Concerning the only not-yet-determined parameter $a_1$, we also discover that it is either normalizable to a real (or imaginary) group parameter or it is never normalizable ({\it cf.} Corollary \ref{cor-a1}). In the former case, the structure group $G$ of the above ambiguity matrices will be reduced to $G^{\sf red}$ of real dimension $1$ while in the later case $G^{\sf red}$ is of real dimension $2$. Next, we start the last part, namely prolongation, of Cartan's method. Accordingly, our equivalence problem to our arbitrary CR model $M_k$ converts by  that to the prolonged space $M_k\times G^{\sf red}$ of real dimension either $3+k$ or $4+k$. Finding the structure equations of this new equivalence problem is easy, it is enough to add the equation $d\alpha=0$ to the above structure equations \thetag{\ref{int-struc-eq}}.

\begin{Theorem}
\label{theorem-introduction}
(cf. Theorem \ref{theorem-main})
The biholomorphic equivalence problem of a totally nondegenerate CR model $M_k$ of codimension $k$ and real dimension $2+k$ is reducible to some absolute parallelisms, namely to some certain $\{e\}$-structures on prolonged manifolds of real dimensions either $3+k$ or $4+k$.
\end{Theorem}

In Section \ref{section-proof}, we start utilizing the achieved results to prove Beloshapka's maximum conjecture \ref{conjecture} in CR dimension one. According to the principles of Cartan's theory (\cite{Olver-1995}), once we receive the final constant type structure equations of the equivalence problem to each CR model $M_k$, we can plainly attain the structure of its symmetry Lie algebra $\frak{aut}_{CR}(M_k)$. Computing this algebra by means of the achieved constant type structure equations \thetag{\ref{int-struc-eq}} together with $d\alpha=0$, we realize that it is graded and without any positive part ({\it cf.} Proposition \ref{prop-graded-algebra}) as was the assertion of the conjecture.

Finally in appendix \ref{appendix}, we illustrate the results by considering the $8$-dimensional length $4$ CR model $M_6$.

It may be worth noting at the end of this section that though we deal mainly in this paper with models of CR dimension one, the results seem  generalizable by means of analogous proofs. Moreover, we find a more general and stronger fact than the assertion of the maximum conjecture about Beloshapka's CR models, namely Theorem \ref{theorem-introduction}. Even more we prove, as a result of inspecting the Lie algebra associated with the above constant type structure equations, that;

\begin{Proposition} (cf. Proposition \ref{prop-graded-algebra})
\label{cor-introduction}
 The Lie algebra $\frak{aut}_{CR}(M_k)$ of a $k$-codimensional weight $\rho$ totally nondegenerate CR model $M_k$ is graded of the form:
 \[
 \frak{aut}_{CR}(M_k):=\underbrace{\frak g_{-\rho}\oplus\ldots\oplus\frak g_{-1}}_{\frak g_-}\oplus\frak g_0
 \]
 where $\frak g_-$ is $(2+k)$-dimensional and where $\frak g_0$ is Abelian of dimension either $1$ or $2$. Thus, we have:
 \[
  {\rm dim}\big(\frak{aut}_{CR}(M_k)\big)= \, 3+k \ \ \ \ or \ \ \ \ \ 4+k.
 \]
\end{Proposition}

 As a homogeneous space, each CR model $M_k$ can be considered as a quotient space ({\it see} the paragraph after equation \thetag{\ref{Aut}}):
\[
M_k\equiv\frac{{\sf Aut}_{CR}(M_k)}{{\sf Aut}_{0}(M_k)}\cong {\sf G_-}
\]
of the CR automorphism group ${\sf Aut}_{CR}(M_k)$, corresponding to $\frak{aut}_{CR}(M_k)$ by its isotropy subalgebra ${\sf Aut}_{0}(M_k)$ at the origin, corresponding to $\frak g_0\oplus\frak g_+$. The above corollary states, in a more precise manner, that such isotropy group is just ${\sf G}_0$, corresponding to the Abelian algebra $\frak g_0$ and comprises only {\it linear} CR automorphisms $h:M_k\rightarrow M_k$, preserving the origin. Even more precisely, in this case that ${\rm dim}\,\frak g_0$ is either $1$ or $2$, then $\sf G_0$ can be identified with the matrix Lie group ${\sf GL}(1,\mathbb R)=(\mathbb R^\ast, \times)$ in the former case and ${\sf GL}(1,\mathbb C)=(\mathbb C^\ast, \times)$ in the latter. This is a decisive generalization of the Theorem, presented at the end of \cite{Mamai2009}.

\section{Beloshapka's models}
\label{section-Beloshapkas-models}

In this preliminary section, we explain  the method of constructing defining equations of Beloshapka's models in CR dimension one. For more  detailed explanation, we refer the reader to \cite{Beloshapka2004, Mamai2013}. In each fixed CR codimension $k$, a certain Beloshapka's model $M_k\subset\mathbb C^{1+k}$ can be represented in coordinates $(z,w_1,\ldots,w_k)$ as the graph of some $k$ real-valued polynomial functions, which {\it roughly speaking}, are the homogeneous parts of the defining equations introduced in Theorem \ref{theorem-existence-totally-momdeg}. Throughout constructing these defining polynomials and to each complex coordinate $\sf x$, it will be assigned a {\sl weight} number $[{\sf x}]$. Recall that for a monomial ${\sf x}_1^{\alpha_1}\cdots {\sf x}_n^{\alpha_n}$, the associated weight is defined as $\sum_{i=1}^n\,\alpha_i\,[{\sf x}_i]$. Moreover, a polynomial is called {\sl weighted homogeneous} of the weight $\sf w$ whenever all of its monomials are of this weight. We assign the same weight $[{\sf x}]$ to the conjugation $\overline{\sf x}$ and real and imaginary parts of $\sf x$, as well. Before starting the construction, we first need the following definition;

\begin{Definition} ({\it cf.} \cite{Krantz})
An arbitrary $\mathcal C^2$ complex function $f:\Omega\subset\mathbb C^n\rightarrow\mathbb C$ in terms of the coordinates $({\sf z}_1,\ldots, {\sf z}_n)$ is called {\sl pluriharmonic} on its domain $\Omega$ whenever for each $i,j=1,\ldots, n$ we have:
\[
\frac{\partial^2 f}{\partial {\sf z}_i\,\partial\overline{{\sf z}}_j}\equiv 0.
\]
\end{Definition}

In the case that $f$ is real-valued, then locally, pluriharmonicity of $f$ is equivalent to state that it is the real part of a holomorphic function (\cite[Propoition 2.2.3]{Krantz}).

 By convention, we assign to the complex variable $z$ the weight $[z]=1$. The weights of the next complex variables $w_1, w_2, \ldots$, which are absolutely bigger than $1$, will be determined as follows, step by step. At the first onset that only the weight of the single variable $z$ is known, let $\mathcal N_2$ be a basis for the space of all non-pluriharmonic real-valued polynomials of the homogeneous weight $2$, in terms of the complex variables $z$ and $\overline z$. A careful inspection shows that $\mathcal N_2$ comprises merely the single term:
\[
\mathcal N_2:=\{z\overline z\}.
\]
Since the cardinality of this set is $k_2=1$, then we assign immediately the weight $2$ to the next one complex variable $w_1$, {\it i.e.} $[w_1]=2$.

At the moment, two of the complex variables $z$ and $w_1$ have received their weight numbers.
Define the next collection $\mathcal N_3$ as a basis for the space of all real-valued polynomials of the weight $3$, in terms of the variables $z,\overline z$ and ${\rm Re}\,w_1$, which are non-pluriharmonic on the submanifold represented by the weight two homogeneous polynomial ${\rm Im}\,w_1=z\overline z$ in $\mathbb C^2$.
Again, a careful inspection shows that:
\[
\mathcal N_3:=\{{\rm Re}\,z^2\overline z=\frac{z^2\overline z+z\overline z^2}{2}, \ \ \ \ {\rm Im}\,z^2\overline z=\frac{z^2\overline z-z\overline z^2}{2i}\}.
\]
This time, since the cardinality of $\mathcal N_3$ is $k_3=2$, then immediately we assign the weight $3$\,\,---\,\,namely the weight of the monomials in $\mathcal N_3$\,\,---\,\,to the next {\it two} complex variables $w_2$ and $w_3$.

Inductively, assume that $\mathcal N_{j_0}$ is the last constructed basis for some integer $j_0\in\mathbb N$. This means that all the complex variables $z, w_1, w_2, w_3, \ldots, w_{\sf r}$ have received their weight numbers where ${\sf r}:=\sum_{i=2}^{j_0}\, k_i$ and where $k_i:={\rm Card}\,\mathcal N_i$. To construct the next collection $\mathcal N_{j_0+1}$ and for the sake of clarity, let us show the $k_i$ elements of each $\mathcal N_i$ as $\mathcal N_i:=\{{\sf t}^i_1, {\sf t}^i_2, \ldots,{\sf t}^i_{k_i}\}$.  Also for each $\ell=2,\ldots, j_0$, let ${\bf w}^\ell=(w_{l}, \ldots, w_{l+k_\ell-1})^t$ be the $k_\ell$-tuple of all complex variables $w_1,\ldots,w_{\sf r}$ of the same weight $\ell$ and consider:
\begin{equation*}
A_\ell=\left(
         \begin{array}{ccc}
           a^\ell_{11} & \ldots & a^\ell_{1k_\ell} \\
           \vdots & \vdots & \vdots \\
            a^\ell_{k_\ell 1} & \ldots & a^\ell_{k_\ell k_\ell} \\
         \end{array}
       \right)
\end{equation*}
 as some real $k_\ell\times k_\ell$ matrix of the maximum ${\rm Rank}(A_\ell)=k_\ell$. Then, the sought collection $\mathcal N_{j_0+1}$ is defined as a basis for the space of all real-valued polynomials of the weight $j_0+1$, in terms of the already weight determined variables $z, \overline z, {\rm Re}\,w_1, {\rm Re}\,w_2, \ldots, {\rm Re}\,w_{\sf r}$, which are non-pluriharmonic on the submanifold represented as the graph of some $\sf r$ weighted homogeneous polynomial functions:
\[
{\rm Im}\,{\bf w}^\ell
=
A_\ell
\cdot
\left(
  \begin{array}{c}
    {\sf t}^\ell_1 \\
    \vdots \\
    {\sf t}^\ell_{k_\ell} \\
  \end{array}
\right) \ \ \ \ \ \ \ \ {\scriptstyle (\ell\,=\,2\,,\,\ldots\,,\,j_0)},
\]
in $\mathbb C^{{\sf r}+1}$. Here, ${\rm Im}\,{\bf w}^\ell$ is the $k_\ell$-tuple of imaginary parts of ${\bf w}^\ell$.
If the cardinality of $\mathcal N_{j_0+1}$ is $k_{j_0+1}$, then one assigns immediately the weight $j_0+1$ to all the next complex variables $w_{{\sf r}+1},\ldots,w_{{\sf r}+k_{j_0+1}}$.

 \subsection{Constructing the defining equations}

After assigning appropriate weights to the complex variables $z, w_\bullet$, we are ready to explain the procedure of constructing  defining polynomials of a $k$-codimensional Beloshapka's model $M_k\subset\mathbb C^{k+1}$. In this case, we need only the assigned weights to the complex coordinates $(z,w_1,\ldots,w_k)$ of $M_k$, hence we have to construct the above sets $\mathcal N_i$ until we arrive at the stage $i=\rho$ where $\rho$ is the smallest integer satisfying:
\begin{equation}
 \label{inequality-k}
 k\leqslant  k_2+\ldots+k_{\rho-1}+k_\rho.
 \end{equation}
 In this case, the chain of associated weights to the complex variables $z,w_1,\ldots,w_k$ is ascending and the last variable $w_k$ is of the maximum weight $\rho$, which will be in fact the length of $M_k$ ({\it cf.} Theorem \ref{theorem-existence-totally-momdeg}).

Now, for each $\ell=2,\ldots, \rho-1$, consider the $k_\ell$-tuple ${\bf w}^\ell$ and the $k_\ell\times k_\ell$ matrix $A_\ell$ as above. For $\ell=\rho$ and since in this case the number of the present weight $\rho$ variables among $w_1,\ldots, w_k$ is $m=k-\sum_{i=2}^{\rho-1}k_i\leqslant k_\rho$, then consider the $m$-tuple ${\bf w}^\rho$ as ${\bf w}^\rho=(w_{k-m+1},\ldots, w_k)$. Also let:
\[
A_\rho=\left(
         \begin{array}{ccc}
           a^\rho_{11} & \ldots & a^\rho_{1,k_\rho} \\
           \vdots & \vdots & \vdots \\
            a^\rho_{m 1} & \ldots & a^\rho_{m k_\rho} \\
         \end{array}
       \right)
\]
be a certain real $m\times k_\rho$ matrix of the maximum ${\rm Rank}(A_\rho)=m$. Then, the desired defining equations of $M_k$ can be represented in the following matrix form:
\begin{equation}
\label{model-matrix}
{\rm Im}\,{\bf w}^\ell
=
A_\ell
\cdot
\left(
  \begin{array}{c}
    {\sf t}^\ell_1 \\
    \vdots \\
    {\sf t}^\ell_{k_\ell} \\
  \end{array}
\right), \ \ \ \ \ \ \ \ {\scriptstyle (\ell\,=\,2\,,\,\ldots\,,\,\rho)}.
\end{equation}

As we observe, in a fixed codimension $k$ one may find infinite number of CR models $M_k$ determined by different values of the above matrix entries $a^\ell_{ij}$. Nevertheless, possibly many of them are equivalent, up to some biholomorphic change of coordinates. For example in codimension $k=3$, CR models $M_3\subset\mathbb C^4$ are represented as the graph of some three defining polynomials:
\[
\aligned
{\rm Im}\,w_1&=a\,z\overline z,
\\
 {\rm Im}\,w_2&= a_{11}\,(z^2\overline z+z\overline z^2)+ia_{12}\,(z^2\overline z-z\overline z^2), \ \ \ \ \ \ \ {\scriptstyle (a,\, \, a_{ij}\,\in\,\mathbb R)},
\\
{\rm Im}\,w_3&= a_{21}\,(z^2\overline z+z\overline z^2)+ia_{22}\,(z^2\overline z-z\overline z^2).
\endaligned
\]
However, by some simple biholomorphic changes of coordinates like those presented at the page 50 of \cite{5-cubic} (expanded version), one shows that they are biholomorphically equivalent to the so-called {\it 5-cubic} model:
\[
\aligned
M^5_{\sf c}: \ \ \ \ \ \
\left[
\begin{array}{l}
{\rm Im}\,w_1=z\overline z,
\\
{\rm Im}\,w_2=z^2\overline z+z\overline z^2,
\\
{\rm Im}\,w_3= i\,(z^2\overline z-z\overline z^2).
\end{array}
\right.
\endaligned
\]

Anyway, in this paper we do not stress on such biholomorphic normalizations since it will not matter whether the under consideration defining equations are normalized or not.

Summing up the above procedure, each arbitrary CR model $M_k\subset\mathbb C^{1+k}$ of codimension $k$ and of the length $\rho$ can be represented as the graph of some $k$ certain real-valued defining equations:
\begin{equation}
\label{model}
\aligned
M_k: \ \ \left\{
\begin{array}{l}
w_1-\overline w_1=2i\,\Phi_1(z,\overline z),
\\
\ \ \ \ \ \ \ \ \ \ \ \ \ \ \ \ \vdots
\\
w_j-\overline w_j=2i\,\Phi_j(z,\overline z,w,\overline w),
\\
\ \ \ \ \ \ \ \ \ \ \ \ \ \ \ \ \vdots
\\
w_k-\overline w_k=2i\,\Phi_k(z,\overline z,w, \overline w),
\end{array}
\right.
\endaligned
\end{equation}
where each $\Phi_j$ is of the weight $[w_j]$, in terms of the complex variables $z,\overline z$ and real variables ${\rm Re}\,w_i$ with $[w_i]\lneqq[w_j]$.
As one observes, the defining equations of $M_k$ are actually those of a certain $(k-1)$-codimensional model $M_{k-1}$, added just by the last equation $w_k-\overline w_k=2i\,\Phi_k(z,\overline z,w,\overline w)$.

\begin{Remark}
Instead of the above Beloshapka's algebraic method for constructing defining equations of a totally nondegenerate CR model $M_k$, Jo\"el Merker in \cite{Merker-5-CR-III} has introduced a more geometric way by considering the affect of the total nondegeneracy on the converging power series expansions of the desired defining equations.
\end{Remark}

\section{Constructing associated structure equations}
\label{section-lifted-coframes}

Studying equivalences between geometric objects by means of Cartan's classical approach entails first some preliminary equipments, the end of them is the construction of associated structure equations. In the current case of biholomorphic equivalence to CR models, we follow the systematic method developed among the recent years in \cite{5-cubic, Samuel, SCM-Moduli, Merker-Sabzevari-M3-C2} which includes three major steps to bring us to the stage of constructing the sought structure equations:

\begin{itemize}
\item[$\bullet$] Finding appropriate initial CR frame for each model and computing its Lie commutators.
\item[$\bullet$] Passage to the dual CR coframe and computing the associated Darboux-Cartan structure.
\item[$\bullet$] Finding the ambiguity matrix of the equivalence problem, in question.
\end{itemize}

\subsection{Associated initial frames for the complexified tangent bundles}
\label{sub-sec-initial-frame}

In the defining equations  \thetag{\ref{model}} of a $k$-codimensional CR model $M_k\subset\mathbb C^{k+1}$ in coordinates $(z,w_1,\ldots,w_k)$, each real-valued polynomial $\Phi_j(z,\overline z,\overline w)$ is ${\rm O}(2)$ and thus we can apply the analytic implicit function theorem in order to solve these equations for the $k$
variables $w_j$, $j = 1, \ldots , k$. Performing this, we obtain equivalently a collection of $k$ complex defining equations like:
\begin{equation}
\label{model-complex}
\aligned
M_k: \ \ \left\{
\begin{array}{l}
w_j=\Theta_j(z,\overline z,\overline w) \ \ \ \ \ \ \ \ \ \ {\scriptstyle (j\,=\,1\,,\,\ldots\,,\, k)},
\end{array}
\right.
\endaligned
\end{equation}
where each {\it complex-valued} polynomial function $\Theta_j$ is in terms of $z, \overline z, \overline w_j$ and some other conjugated variables $\overline w_\bullet$ of absolutely lower weights than $[w_j]$. By an induction on the weights associated with the complex coordinates $w_1, \ldots, w_k$, one verifies that similar to the real-valued functions $\Phi_\bullet$ also each complex-valued polynomial $\Theta_j$ is weighted homogeneous of the weight $[w_j]$.

Having in hand the complex defining polynomials \thetag{\ref{model-complex}} of the CR model $M_k$ and according to \cite{Merker-Porten-2006, 5-cubic}, then the associated holomorphic and antiholomorphic tangent bundles $T^{1,0}M_k$ and $T^{0,1}M_k$, can be generated respectively by the single vector fields:

\begin{equation}
\label{L1}
\aligned
\mathcal L:=\frac{\partial}{\partial z}+\sum_{j=1}^k\,\frac{\partial \Theta_j}{\partial z}(z,\overline z, \overline w)\,\frac{\partial}{\partial w_j} \ \ \ \ \textrm{and} \ \ \ \ \overline{\mathcal L}:=\frac{\partial}{\partial \overline z}+\sum_{j=1}^k\,\frac{\partial \overline\Theta_j}{\partial \overline z}(z,\overline z, w)\,\frac{\partial}{\partial\overline w_j}.
\endaligned
\end{equation}

 For ${\sf x}$ to be one of the complex variables $z,w_1,\ldots,w_k$, or one of their conjugations, or one of their real or imaginary parts, we assign the weight $-[{\sf x}]$ to the standard vector filed $\frac{\partial}{\partial {\sf x}}$.
Notice that for a weighted homogeneous polynomial $F({\sf x}, \overline{\sf x})$, each differentiation of the shape $F_{{\sf x}_i}$ (or $F_{\overline{\sf x}_i}$) decreases its weight by $[{\sf x}_i]$ numbers, if it does not vanish. Then, by a glance on the above expressions of $\mathcal L$ and $\overline{\mathcal L}$ one finds them as two weighted homogeneous fields of the same weight $-1$.

 \subsubsection{Notations}

Henceforth and in order to stress their lengths (and weights), let us denote by $\mathcal L_{1,1}$ and $\mathcal L_{1,2}$ the above vector field $\mathcal L$ and its conjugation $\overline {\mathcal L}$, respectively. By the total nondegeneracy of our length $\rho$  fixed CR model $M_k$, one constructs the sought initial frame on $\mathbb C\otimes TM_k$ by applying the iterated Lie brackets\,\,---\,\,or simple words in the terminology of free Lie algebras\,\,---\,\,of these two vector fields, up to the length $\rho$. Let us denote by $\mathcal L_{\ell,i}$ and call it by the $i$-th {\it initial vector field}, the $i$-th appearing independent vector filed obtained as an iterated Lie bracket of the length $\ell$. For example, the next and third appearing vector filed can be computed as the length two iterated Lie bracket:
\[
\mathcal L_{2,3}=[\mathcal L_{1,1},\mathcal L_{1,2}].
\]
In the case that the reference to the order $i$ of a length $\ell$ initial vector field $\mathcal L_{\ell,i}$ is superfluous and by abuse of notation, we denote it just by $\mathcal L_\ell$, which actually is a vector field expressible (inductively) as:
 \begin{equation}
 \label{iterated-bracket}
\mathcal L_\ell:=[\mathcal L_{1,i_1},\underbrace{[\mathcal L_{1,i_2},[\ldots,[\mathcal L_{1,i_{\ell-1}},\mathcal L_{1,i_\ell}]]]}_{\mathcal L_{\ell-1}}] \ \ \ \ \ \ \scriptstyle{(i_j=1,2, \ \ \ell=1,\ldots,\rho)}.
 \end{equation}
 Notice that according to the expressions of $\mathcal L_{1,1}$ and $\mathcal L_{1,2}$ in  \thetag{\ref{L1}} and from the length $\ell= 2$ to the end, one does not see any coefficient of $\frac{\partial}{\partial z}$ or $\frac{\partial}{\partial\overline z}$ in the expression of $\mathcal L_\ell$.

\begin{Lemma}
\label{homogeneous}
Each length $\ell$ initial vector field $\mathcal L_{\ell}$ is homogeneous of the weight $-\ell$ with polynomial coefficients. Moreover, for  two initial vector fields  $\mathcal L_{\alpha,i}$ and $\mathcal L_{\beta,j}$ with $\alpha+\beta=\ell$, if $[\mathcal L_{\alpha,i}, \mathcal L_{\beta,j}]\not\equiv 0$ then it is a weighted homogeneous vector field of the weight $-\ell$, again with polynomial coefficients.
\end{Lemma}

\proof
Since the coefficients of the basis vector fields $\mathcal L_{1}$ are of polynomial type, the polynomiality of the coefficients in their iterated brackets is obvious.  Concerning the weights, we continue by a plain induction on the length $\ell$. As we saw, the two vector fields $\mathcal L_{1,1}$ and $\mathcal L_{1,2}$ of the length $\ell=1$ are of the homogeneous weight $-1$. For the next lengths and as our induction hypothesis, assume that all length $\ell$ vector fields:
\[
\mathcal L_{\ell}:=\sum_{[w_i]\geq\ell}\,\varphi_i(z,\overline z, w,\overline w)\,\frac{\partial}{\partial w_i}+\sum_{[w_i]\geq\ell}\,\psi_i(z,\overline z, w,\overline w)\,\frac{\partial}{\partial\overline w_i}
\]
 are weighted homogeneous of the weight $-\ell$. Thus, the nonzero polynomial coefficients $\varphi_i$ and $\psi_i$ are homogeneous of the nonnegative weights $[w_i]-\ell$. Now, consider an arbitrary new appearing initial field $\mathcal L_{\ell+1}=[\mathcal L_{1},\mathcal L_\ell]$ of the length $\ell+1$. Applying the Leibniz rule on this bracket with the present expressions of the weighted homogeneous fields $\mathcal L_1$ in \thetag{\ref{L1}} and $\mathcal L_\ell$, it manifests itself as a weight $-(\ell+1)$ homogeneous vector field, as was expected. The proof of the second part of the assertion is completely similar.
 \endproof

 \begin{Proposition}
\label{prop-neessary-for-Darboux}
Let $\mathcal L_{\alpha,i}$ and $\mathcal L_{\beta,j}$ be two initial vector fields associated with a length $\rho$ CR model $M_k$ with $\alpha+\beta=\ell$. Then, we have:
\[
[\mathcal L_{\alpha,i}, \mathcal L_{\beta,j}]=\sum_{t}\, {\sf c}_{t}\,\mathcal L_{\ell,t},
\]
for some constant integers ${\sf c}_t$. In particular if $\ell>\rho$, then $[\mathcal L_{\alpha,i}, \mathcal L_{\beta,j}]\equiv 0$.
\end{Proposition}

\proof
Since our length $\rho$ CR model $M_k$ is totally nondegenerate and according to the discussion after Definition \ref{totally-nondegenerate}, one verifies that the equality holds in the case that $\ell=1,\ldots,\rho-1$. Hence, let us prove the assertion first by assuming $\ell=\rho$. If $[\mathcal L_{\alpha,i}, \mathcal L_{\beta,j}]$ vanishes, then it remains nothing to prove. Otherwise, Lemma \ref{homogeneous} indicates that this Lie bracket produces a weighted homogeneous vector field of the weight $-\rho$ with polynomial coefficients. Taking into account that the maximum weight of the extant complex variables is $\rho$, then this bracket can be written just in the form:
 \[
 [\mathcal L_{\alpha,i}, \mathcal L_{\beta,j}]:=\sum_{[w_i]=\rho}\,{\sf a}_i\,\frac{\partial}{\partial w_i} + \sum_{[w_i]=\rho}\,{\sf b}_i\,\frac{\partial}{\partial \overline w_i},
 \]
 for some constant integers ${\sf a}_i$ and ${\sf b}_i$.
 On the other hand, the projection map $\pi:\mathbb C^{1+k}\rightarrow\mathbb R^{2+k}$ defined as $(x,y,{\bf u}, {\bf v})\mapsto (x,y,{\bf u})$ constitutes a natural local chart-map on the real submanifold $M_k$\,\,---\,\,still we denote $z=x+iy$ and $w_j=u_j+i\,v_j$. Then, one restates {\it intrinsically} the expression of each $\mathcal L_{\ell}$ in terms of $z, \overline z, u_j$ by dropping $\frac{\partial}{\partial v_j}$ for $j=1,\ldots, k$ and also replacing each $v_j$ by its expression in \thetag{\ref{model}}. Then in particular, each weight $-\rho$ initial vector field $\mathcal L_{\rho}$ is expressible as some combinations, with constant coefficients, of standard fields $\frac{\partial}{\partial u_j}$ with $[u_j]=[w_j]=\rho$. Similar fact holds also for the above bracket $[\mathcal L_{\alpha,i}, \mathcal L_{\beta,j}]$. According to Theorem \ref{theorem-existence-totally-momdeg}, the number of standard fields $\frac{\partial}{\partial u_j}$, with $[u_j]=\rho$ is exactly equal to the number of (linearly independent) initial vector fields of the length $\rho$. This implies that the $\mathbb C$-vector space generated by all standard fields $\frac{\partial}{\partial u_j}$ with $[u_j]=\rho$ is equal to that generated by the length $\rho$ initial vector fields $\mathcal L_{\rho}$. Consequently, as an element of this space, the Lie bracket $[\mathcal L_{\alpha,i}, \mathcal L_{\beta,j}]$ can be expressed as a linear combination of the length $\rho$ initial vector fields $\mathcal L_\rho$ with constant coefficients.

  \noindent
 To continue the proof, now let $\ell=\alpha+\beta>\rho$ and suppose, to derive a contradiction, that $[\mathcal L_{\alpha,i}, \mathcal L_{\beta,j}]\not\equiv 0$. Then, according to Lemma \ref{homogeneous}, it should be a weight $-\ell$ vector field with polynomial coefficients:
 \begin{equation*}
[\mathcal L_{\alpha,i}, \mathcal L_{\beta,j}]:=\sum_{i}\,\varphi_i(z,\overline z, w, \overline w)\,\frac{\partial}{\partial w_i}+\sum_{i}\,\psi_i(z,\overline z, w, \overline w)\,\frac{\partial}{\partial\overline w_i},
\end{equation*}
where, consequently, each polynomial $\varphi_i$ and $\psi_i$ is of the weight $[w_i]-\ell$. Since these coefficients are of the polynomial type, then their weights are nonnegative and thus $[w_i]\geq\ell>\rho$. This is a contradiction to the fact that among our coordinates there is no any complex variable of the weight absolutely bigger than $\rho$.
\endproof

\begin{Remark}
According to the notations introduced before Definition \ref{totally-nondegenerate}, the complexified tangent bundle $\mathbb C\otimes TM_k$ of an arbitrary length $\rho$ CR model $M_k$ of codimension $k$, admits the filtration:
 \[
 T^{1,0}M_k+T^{0,1}M_k=D_1\subset D_2\subset\ldots\subset D_\rho=\mathbb C\otimes TM_k,
 \]
 where each $D_\ell$ is a subdistribution constituted by initial vector fields $\mathcal L_l$ of the lengths $l\leqslant\ell$. The above proposition indicates that at each point $p\in M_k$ near the origin, $\mathbb C\otimes T_pM_k$ is isomorphic to the {\it graded complex nilpotent Lie algebra}:
 \[
\frak m:=\frak m_{-\rho}\oplus\frak m_{-\rho+1}\oplus\ldots\oplus\frak m_{-1}
 \]
 where $\frak m_{-1}:=D_1$ and where $\frak m_{-\ell}:=D_\ell\setminus D_{\ell-1}, \ell=2,\ldots,\rho$ is the $\mathbb C$-vector space generated by all initial vector fields of the precise weight $-\ell$. In this case, $D_1$ is called a distribution {\sl of the constant type} $\frak m$.
\end{Remark}

\subsection{The associated initial coframes and their Darboux-Cartan structures}
\label{Darboux section}

For $\ell=1,\ldots,\rho$ and $i=1,\ldots,2+k$, let us denote by $\sigma_{\ell,i}$  the dual CR 1-form associated with the initial vector field $\mathcal L_{\ell,i}$. Since the collection of the weighted homogeneous vector fields $\{\mathcal L_{1,1}, \ldots,\mathcal L_{\rho,2+k}\}$ forms a frame for the complexified bundle $\mathbb C\otimes TM_k$, then its dual set $\{\sigma_{1,1}, \ldots,\sigma_{\rho,2+k}\}$ is a coframe for it.

\begin{Lemma}
\label{lem-d}
Given a frame $\big\{ \mathcal{ V}_1, \dots, \mathcal{
V}_n\big\}$ on an open subset of $\mathbb R^n$ enjoying the Lie structure:
\[
\big[\mathcal{V}_{i_1},\,\mathcal{V}_{i_2}\big] =
\sum_{k=1}^n\,c_{i_1,i_2}^k\,\mathcal{V}_k \ \ \ \ \ \ \ \ \ \ \ \ \
{\scriptstyle{(1\,\leqslant\,i_1\,<\,i_2\,\leqslant\,n)}},
\]
where the $c_{ i_1, i_2}^k$ are certain functions on $\mathbb R^n$, the dual
coframe $\{ \omega^1, \dots, \omega^n \}$ satisfying by definition:
\[
\omega^k
\big(\mathcal{V}_i\big)
=
\delta_i^k
\]
enjoys a quite similar Darboux-Cartan structure, up to an overall minus sign:
\[
d\omega^k = - \sum_{1\leqslant i_1<i_2\leqslant n}\,
c_{i_1,i_2}^k\,\omega^{i_1}\wedge\omega^{i_2}
\ \ \ \ \ \ \ \ \ \ \ \ \ {\scriptstyle{(k\,=\,1\,\cdots\,n)}}.
\]
\end{Lemma}

As a direct consequence of the above Lemma and Proposition \ref{prop-neessary-for-Darboux}, we find the Darboux-Cartan structure of our initial coframe as follows;

\begin{Proposition}
\label{Darboux-Cartan}
The  exterior differentiation of each $1$-form $\sigma_\ell$ dual to the weight $-\ell$ initial vector field $\mathcal L_\ell$ is of the form:
\[
d\sigma_\ell:=\sum_{\beta+\gamma=\ell}\,{\sf c}_{\beta,\gamma}\,\sigma_\beta\wedge\sigma_\gamma,
\]
for some constant complex integers ${\sf c}_{\beta,\gamma}$. This equivalently means that in the expression of each corresponding Lie bracket $[\mathcal L_\beta,\mathcal L_\gamma]$, with $\beta+\gamma=\ell$, the coefficient of $\mathcal L_\ell$ is $-{\sf c}_{\beta,\gamma}$.
\end{Proposition}

\medskip
\noindent
{\bf Weight assignment.} Naturally, we assign the weight $-\ell$ to a certain 1-form $\sigma_{\ell,i}$ and its differentiation $d\sigma_{\ell,i}$ as is the weight of their corresponding field $\mathcal L_{\ell,i}$. Also, we occasionally say that $\sigma_{\ell,i}$ is of the length $\ell$.

Another simple but quite useful result is as follows;

\begin{Lemma}
\label{lem-unique-1-form}
For each weight $-\ell$ initial $1$-form $\sigma_{\ell,i}$ with $\ell\neq1$, there is a weight $-(\ell-1)$ initial $1$-form $\sigma_{\ell-1,j}$ where either $\sigma_{\ell-1,j}\wedge\sigma_{1,1}$ or $\sigma_{\ell-1,j}\wedge\sigma_{1,2}$ is visible {\sl uniquely} in the Darboux-Cartan structure of $d\sigma_{\ell,i}$.
\end{Lemma}

\proof
This is a straightforward consequence of the fact that in the procedure of constructing our initial frame, each weight $-\ell$ vector field $\mathcal L_{\ell,i}$ is constructed as the Lie bracket between $\mathcal L_{1,1}$ or $\mathcal L_{1,2}$ and a {\it unique} weight $-(\ell-1)$ vector filed $\mathcal L_{\ell-1,j}$. Then, Lemma \ref{lem-d} implies the desired results.
\endproof

\subsection{Ambiguity matrix}
\label{sec-ambiguity-matrix}

After providing the above appropriate initial frame and coframe on the complexified tangent bundle $\mathbb C\otimes TM_k$, now this is the time of seeking the ambiguity matrix associated with the problem, what actually encodes biholomorphic equivalences to $M_k$. The procedure of construction is demonstrated in the recent works \cite{Merker-Sabzevari-M3-C2, Samuel, 5-cubic, SCM-Moduli} in the specific cases of $k=1,2,3,4$.
Let us explain it here in the general case of the CR model $M_k$. Assume that:
\[
\aligned
h\colon \ \ \ M_k&\longrightarrow {\bf M}_k
\\
(z,w)
& \longmapsto
\big(z'(z,w),\,w'(z,w)\big)
\endaligned
\]
is a (biholomorphic) equivalence map between our $(2+k)$-dimensional CR model $M_k$
and another arbitrary real analytic totally nondegenerate CR generic submanifold ${\bf M}_k \subset \mathbb{C}^{1+k}$ of codimension $k$, in canonical coordinates $\big( z', w_1', \ldots,
w_k'\big)$.
We assume that
${\bf M}_k$ is also equipped with a frame of $2+k$ {\it lifted} vector fields $\{{\bf L}_{1,1},\,{\bf L}_{1,2}, \, {\bf L}_{2,3}, \, {\bf L}_{3,4}, \, {\bf L}_{3,5}, \ldots, {\bf L}_{\rho,2+k}\}$ where, as before, ${\bf L}_{1,1}$ and ${\bf L}_{1,2}=\overline {{\bf L}_{1,1}}$ are local generators of $T^{1,0} {\bf M}_k$ and $T^{0,1}{\bf M}_k$ and where each other vector field ${\bf L}_{\ell,i}$ can be computed as an iterated Lie bracket between ${\bf L}_{1,1}$ and ${\bf L}_{1,2}$ of the length $\ell$, exactly as \thetag{\ref{iterated-bracket}} for constructing the initial vector filed $\mathcal L_{\ell,i}$. Tensoring with $\mathbb C$, then the push-forward $h_\ast
\colon\ \ \
TM_k
\longrightarrow
T{\bf M}_k$ of $h$ induces a complexified map, still
denoted by the same symbol with the customary abuse of notation (\cite{Boggess-1991}):
\[
\aligned
h_\ast
\colon\ \ \
\mathbb{C}&\otimes
TM_k\longrightarrow
\mathbb{C}
\otimes
T{\bf M}_k,
\\
{\sf z}&\otimes \mathcal X
\longmapsto
{\sf z}\otimes
h_\ast(\mathcal X).
\endaligned
\]
Our current purpose is to seek the associated matrix to this linear map.

According to principles in CR geometry (\cite{BER, Boggess-1991, Merker-Pocchiola-Sabzevari-5-CR-II}), $h_\ast$ transfers every generator of $T^{1,0}M_k$ to a vector field in the same bundle $T^{1,0}{\bf M}_k$. Hence for the single generator $\mathcal L_{1,1}$ of $T^{1,0}M_k$, there exists some {\it nonzero} function $a_1:=a_1(z',w')$ with:
\begin{equation}
\label{h*-1}
h_\ast(\mathcal L_{1,1})=a_1\,{\bf L}_{1,1}.
\end{equation}
Moreover, $h_\ast$ preserves the conjugation, whence for $\mathcal L_{1,2}:=\overline{\mathcal L_{1,1}}$, we have $h_\ast(\mathcal L_{1,2})=\overline a_1\,{\bf L}_{1,2}$.

The third vector field in the basis of $\mathbb C\otimes TM_k$ is the {\it imaginary field} $\mathcal L_{2,1}:=[\mathcal L_{1,1}, \mathcal L_{1,2}]$. Let us compute the image of $h_\ast$ on it:
\begin{equation}
\label{h*-2}
\aligned
h_{\ast}(\mathcal L_{2,3})
&= h_{\ast}\big([\mathcal L_{1,1},\mathcal L_{1,2}]\big)=
\big[h_\ast(\mathcal L_{1,1}),h_\ast(\mathcal L_{1,2})\big]
=\big[a_1\,{\bf L}_{1,1},\overline a_1\,{\bf L}_{1,2}\big]
\\
&
=
a_1\overline a_1 \big[{\bf L}_{1,1}\,{\bf L}_{1,2}\big]
\underbrace{-
\overline a_1\,{\bf L}_{1,2}(a_1)}_{
=:\,a_2}
{\bf L}_{1,1}
+
a_1\,{\bf L}_{1,1}\big(\overline a_1\big)
{\bf L}_{1,2}
\\
&=: a_1\overline a_1{\bf L}_{2,3}+a_2\,{\bf L}_{1,1}-\overline a_2\,{\bf L}_{1,2},
\endaligned
\end{equation}
for a certain function $a_2:=a_2(z',w')$.

Next, in the length three, two initial fields $\mathcal L_{3,4}=[\mathcal L_{1,1},\mathcal L_{2,3}]$ and $\mathcal L_{3,5}=-[\mathcal L_{1,2},\mathcal L_{2,3}]$ exist. Without affecting the results, we multiply the second bracket by $-1$ to have the simple relation $\mathcal L_{3,5}=\overline {\mathcal L_{3,4}}$. In a similar fashion of computations, one finds:

\begin{equation}
\label{h*-3}
\aligned
h_\ast(\mathcal L_{3,4})&:=
a_1^2\overline a_1\,{\bf L}_{3,4}+\underbrace{\big(a_1\,{\bf L}_{1,1}\big(a_1\overline a_1\big)-a_1\overline a_2\big)}_{=:a_3}{\bf L}_{2,3}+
\\
&+\underbrace{\big(-{a_1\overline a_1}\,{\bf L}_{2,3}(a_1)+ a_1\,{\bf L}_{1,1}(a_2)-a_2\,{\bf L}_{1,1}(a_1)+\overline a_2\,{\bf L}_{1,2}(a_1)\big)}_{=:a_4}{\bf L}_{1,1}
\underbrace{-a_1{\bf L}_{1,1}\big(\overline a_2\big)}_{=:a_5}{\bf L}_{1,2},
\endaligned
\end{equation}
for some three certain complex functions $a_j:=a_j(z',w'), \ j=3,4,5$. By conjugation, we also have:
\[
h_\ast(\mathcal L_{3,5})=a_1\overline a_1^2\,{\bf L}_{3,5}-\overline a_3\,{\bf L}_{2,3}+\overline a_5\,{\bf L}_{1,1}+\overline a_4\,{\bf L}_{1,2}.
\]

Proceeding along the same lines of computations and by an induction on the weight of the initial fields, one finds the following general expression for the image of the complexified push-forward map $h_\ast$;

\begin{Lemma}
\label{lem-number-L11-12}
For a fixed length $\ell$ initial vector field ${\mathcal L}_{\ell,i}$, the push-forward map $h_\ast$ transfers it to a combination like:
\begin{equation}
\label{h*-general}
h_\ast({\mathcal L}_{\ell,i}):=a_1^{p}\overline a_1^q\,{\bf L}_{\ell,i}+\sum_{l<\ell}\,{\sf a}_{r_j}\,{\bf L}_{l,r}, \ \ \ \ \ {\rm with} \ \ \ p+q=\ell
\end{equation}
where ${\sf a}_{r_j}$s are some (possibly zero) complex functions in terms of the target coordinates $(z',w')$.
In other words, $h_\ast({\mathcal L}_{\ell,i})$ is a combination of the corresponding lifted vector field  ${\bf L}_{\ell,i}$ and some other ones ${\bf L}_{l,r}$ of absolutely smaller lengths $l<\ell$.
\end{Lemma}

Then, our sought invertible matrix associated with $h_\ast$ is a $(2+k)\times (2+k)$ upper triangular matrix satisfying\,\,---\,\,here we drop the push-forward $h_\ast$ at the left hand side, for simplicity:

\begin{equation}
\label{ambiguity-frame-k}
\footnotesize\aligned
\left(
  \begin{array}{c}
    \mathcal L_{\rho,i} \\
    \mathcal L_{\rho-1, j} \\
    \vdots \\
    \mathcal L_{3,5} \\
    \mathcal L_{3,4}\\
    \mathcal L_{2,3} \\
    \mathcal L_{1,2} \\
    \mathcal L_{1,1} \\
  \end{array}
\right)
=
\left(
  \begin{array}{cccccccc}
    a_1^{p}\overline a_1^{q} & {\sf a}_\bullet & {\sf a}_\bullet & {\sf a}_\bullet & {\sf a}_\bullet & {\sf a}_\bullet & {\sf a}_\bullet & {\sf a}_\bullet \\
    0 & a_1^{p'}\overline a_1^{q'} & {\sf a}_\bullet & {\sf a}_\bullet & {\sf a}_\bullet & {\sf a}_\bullet & {\sf a}_\bullet & {\sf a}_\bullet \\
    0 & 0 & \ddots & {\sf a}_\bullet & \ldots & \ldots & \ldots & {\sf a}_\bullet \\
    0 & \ldots & 0 & a_1\overline a_1^2 & 0 & -\overline a_3 & \overline a_4 & \overline a_5 \\
    0 & \ldots & 0 & 0 & a_1^2\overline a_1 & a_3 & a_5 & a_4 \\
    0 & 0 & \ldots & \ldots & 0  & a_1\overline a_1 & -\overline a_2 & a_2 \\
    0 & 0 & 0 & \ldots & \ldots & 0 & \overline a_1 & 0 \\
    0 & 0 & 0 & 0 & \ldots & \ldots & 0 & a_1 \\
  \end{array}
\right)
\cdot
\left(
  \begin{array}{c}
    {\bf L}_{\rho,i} \\
    {\bf L}_{\rho-1,j} \\
    \vdots \\
    {\bf L}_{3,5} \\
    {\bf L}_{3,4}\\
    {\bf L}_{2,3} \\
    {\bf L}_{1,2} \\
    {\bf L}_{1,1} \\
  \end{array}
\right).
\endaligned
\end{equation}
If on the main diagonal of the matrix and in front of $\mathcal L_{\ell,r}$ we have $a^{p_r}\overline a^{q_r}$, then $p_r+q_r=\ell$.
As a result of {\it explicitness} in the already procedure of constructing the above desired matrix, we have also the following key observation;

\begin{Lemma}
\label{lem-a2-a3}
In the case that both the specific group parameters $a_2$ and $a_3$, appeared in \thetag{\ref{h*-2}} and \thetag{\ref{h*-3}}, vanish then all the next parameters $a_4, a_5, \ldots$ vanish, identically.
\end{Lemma}
\proof
First we claim that all the appearing group parameters $a_j$ with $j>1$ are some combinations of the iterated $\{{\bf L}_{1,1}, {\bf L}_{1,2}\}$-differentiations of the first parameter $a_1$ and its conjugation $\overline a_1$. We prove our claim by an induction on the length of the initial fields. By \thetag{\ref{h*-2}}, the claim holds for $a_2$ and as our induction hypothesis, assume that it holds for all group parameters $a_\bullet$ appearing among computing the image of $h_\ast$ on initial fields $\mathcal L_l$ of the lengths $\leqslant\ell$. Then for an initial field $\mathcal L_{\ell+1}=[\mathcal L_{1,1},\mathcal L_{\ell,i}]$ of the next length $\ell+1$ (similar argument holds if we have $\mathcal L_{1,2}$ in place of $\mathcal L_{1,1}$) and according to the above Lemma \ref{lem-number-L11-12} we have:
 \[
 h_\ast(\mathcal L_{\ell+1})=\big[h_\ast(\mathcal L_{1,1}),h_\ast(\mathcal L_{\ell,i})\big]=\big[a_1{\bf L}_{1,1},a_1^{p}\overline a_1^q\,{\bf L}_{\ell,i}+\sum_{l<\ell}\,{\sf a}_{r_j}\,{\bf L}_{l,r}\big],
 \]
 where, by hypothesis induction, the appearing coefficients ${\sf a}_{r_j}$ are some combinations of the iterated $\{{\bf L}_{1,1}, {\bf L}_{1,2}\}$-differentiations of $a_1$ and $\overline a_1$. Computing this bracket by means of the Leibniz rule, one finds the new coefficients, namely new group parameters, again as some combinations of the iterated $\{{\bf L}_{1,1}, {\bf L}_{1,2}\}$-differentiations of $a_1$ and $\overline a_1$. This completes the proof of our claim.

 \noindent
 Now, according to \thetag{\ref{h*-2}} and \thetag{\ref{h*-3}} we have:
\[
\aligned
a_2=-\overline a_1\,{\bf L}_{1,2}(a_1) \ \ \ \ {\rm and} \ \ \ \ a_3=a_1\,{\bf L}_{1,1}\big(a_1\overline a_1\big)-a_1\overline a_2.
\endaligned
\]
Since $a_1\neq 0$, then vanishing of $a_2$ and $a_3$ implies that\,\,---\,\,reminding ${\bf L}_{1,2}=\overline {{\bf L}_{1,1}}$:
\[
  {\bf L}_{1,1}(\overline a_1)\equiv 0, \ \ \ \ {\bf L}_{1,1}(a_1)\equiv 0, \ \ \ \ {\bf L}_{1,2}(a_1)\equiv 0, \ \ \ \ {\bf L}_{1,2}(\overline a_1)\equiv 0.
\]
Thus according to our claim, if $a_2$ and $a_3$ vanish then all the next group parameters $a_j$ vanish, identically.
\endproof

\medskip
\noindent
{\bf Weight assignment.} Let $a_j$ be a group parameter which is appeared among computing the value of $h_\ast$ on a length $\ell$ initial vector field $\mathcal L_\ell$. Then, we assign the weight $\ell$ to this group parameter and its conjugation $\overline a_j$. For example, according to \thetag{\ref{h*-1}}, \thetag{\ref{h*-2}} and \thetag{\ref{h*-3}} we have:
\[
[a_1]=1, \ \ \ [a_2]=2, \ \ \ [a_3]=[a_4]=[a_5]=3.
\]
  By this assignment, the nonzero entries at each row of the above matrix \thetag{\ref{ambiguity-frame-k}} have equal weight.

 For each lifted vector field ${\bf L}_{\ell,i}$, let us denote by $\Gamma_{\ell,i}$ its dual {\it lifted} $1$-form and as its corresponding initial 1-form $\sigma_{\ell,i}$, assign the weight $-\ell$ to it. The sough ambiguity matrix $\bf g$ of our equivalence problem, in question, is defined as the invertible matrix associated with the dual pull-back $h^\ast:\mathbb C\otimes T^\ast{\bf M}_k\rightarrow\mathbb C\otimes T^\ast M_k$ of the push-forward $h_\ast$. Then, after a plain matrix transposition we have:

\begin{equation}
\label{ambiguity-coframe-k}
\footnotesize\aligned
\left(
  \begin{array}{c}
    \Gamma_{\rho,i}\\
    \vdots \\
    \Gamma_{\rho-1,j} \\
    \vdots \\
    \Gamma_{3,5} \\
    \Gamma_{3,4} \\
    \Gamma_{2,3} \\
    \Gamma_{1,2} \\
    \Gamma_{1,1} \\
  \end{array}
\right)
=
\underbrace{\left(
  \begin{array}{cccccccc}
    a_1^{p}\overline a_1^{q} & 0 & 0 & 0 & 0 & 0 & 0 & 0 \\
    {\bf 0} & \vdots & \vdots & \vdots & \vdots & \vdots & \vdots & \vdots \\
    {\sf a}_\bullet & a_1^{p'}\overline a_1^{q'} & 0 & 0 & \ldots & \ldots & 0 & 0 \\
    {\sf a}_\bullet & {\sf a}_\bullet & \ddots & 0 & 0 & \ldots & \ldots & 0 \\
    {\sf a}_\bullet & {\sf a}_\bullet & \ldots & a_1\overline a_1^2 & 0 & 0 & \ldots & 0 \\
    {\sf a}_\bullet &\ldots & \ldots & 0 & a_1^2\overline a_1 & 0 & \ldots & 0 \\
    {\sf a}_\bullet & \ldots & \ldots & -\overline a_3 & a_3 & a_1\overline a_1 & 0 & 0 \\
    {\sf a}_\bullet & \ldots & \ldots & \overline a_4 & a_5 & -\overline a_2 & \overline a_1 & 0 \\
    {\sf a}_\bullet & {\sf a}_\bullet & \ldots & \overline a_5 & a_4 & a_2 & 0 & a_1 \\
  \end{array}
\right)}_{{\bf g}}
\cdot
\left(
  \begin{array}{c}
    \sigma_{\rho,i} \\
    \vdots \\
    \sigma_{\rho-1,j} \\
    \vdots \\
    \sigma_{3,5} \\
    \sigma_{3,4} \\
    \sigma_{2,3} \\
    \sigma_{1,2} \\
    \sigma_{1,1} \\
  \end{array}
\right).
\endaligned
\end{equation}

\begin{Remark}
Clarifying the structure of the matrix $\bf g$, it is important to notice that thanks to Lemma \ref{lem-number-L11-12} and for each arbitrary $i$-th {\it column} of this matrix, the first nonzero entry, which stands at the diagonal, is of the form $a_1^r\overline a_1^s$. Even more, since the only length $\ell$ lifted vector field in the image $h_\ast(\mathcal L_{\ell,i})$ in \thetag{\ref{h*-general}} is ${\bf L}_{\ell,i}$, we can state that: if the $i$-th {\it row} of the left (or right) hand side vertical matrix in \thetag{\ref{ambiguity-coframe-k}} is of the weight $-\ell$, then all the entries at the $i$-th column of $\bf g$ standing below $a_1^r\overline a_1^s$ and in front of a weight $-\ell$ 1-form $\Gamma_\ell$ are zero. This fact is shown for example by the zero vector $\bf 0$ in the first column of $\bf g$  or by the entry $0$ below $a_1\overline a_1^2$.
\end{Remark}

\begin{Lemma}
\label{lemma-weight}
If the $1$-form at the $i$-th {\sl row} of the left (or right) hand side vertical matrix of \thetag{\ref{ambiguity-coframe-k}} is of the weight $-\ell$ then, all nonzero entries at the $i$-th {\sl column} of $\bf g$ are of the same weight $\ell$.
\end{Lemma}

\proof
It is a straightforward consequence of the two paragraphs mentioned before \thetag{\ref{ambiguity-coframe-k}}.
\endproof

The collection of all invertible matrices of the form $\bf g$ constitutes a finite dimensional (matrix) Lie group $G$, called by the {\it structure Lie group} of the equivalence problem to the CR model $M_k$.

Recall that ({\it see} \thetag{\ref{model}} and the paragraph after it) the defining equations of our $k$-codimensional CR model $M_{k}\subset\mathbb C^{1+k}$ are precisely those of a CR model  $M_{k-1}$ of codimension $k-1$  added just by the last equation $w_k-\overline w_k=2i\Phi_k(z,\overline z, w, \overline w)$. Let us state a result that will be of use later;

\begin{Proposition}
\label{prop-g-k-1}
The $(1+k)\times (1+k)$ ambiguity matrix ${\bf g}_{k-1}$ associated with the CR model $M_{k-1}$ is a submatrix of the ambiguity matrix $\bf g$ associated with $M_k$, standing as ({\it cf.} \thetag{\ref{ambiguity-coframe-k}}):
 \begin{equation}
 \label{g-k-1}
 \footnotesize\aligned
 {\bf g}=\left(
  \begin{array}{cc}
    a_1^{p}\overline a_1^{q} & 0 \,\,\,\, 0 \,\,\,\, \cdots \,\,\,\, 0 \,\,\,\, 0 \,\,\,\, 0 \,\,\,\, 0  \\
    \begin{array}{c}
      {\bf 0} \\
      {\sf a}_\bullet \\
      \vdots \\
      {\sf a}_\bullet \\
      {\sf a}_\bullet \\
      {\sf a}_\bullet \\
      {\sf a}_\bullet
    \end{array}
         &
     \left(
       \begin{array}{ccccccc}
          &  &  &  &  &  &  \\
          &  &  &  &  &  &  \\
          &  &  &  &  &  &  \\
          &  &  & {\bf g}_{k-1} &  &  &  \\
          &  &  &  &  &  &  \\
          &  &  &  &  &  &  \\
          &  &  &  &  &  &  \\
       \end{array}
     \right)
           \\
  \end{array}
\right).
 \endaligned
 \end{equation}
\end{Proposition}
\proof
Let $M_{k-1}$ be of the length $\rho'\leq\rho$. If we proceed {\it ab initio} as subsection \ref{sub-sec-initial-frame} to provide an initial frame $\{\mathcal L^{\sf old}_{1,1}, \mathcal L^{\sf old}_{1,2}, \ldots, \mathcal L^{\sf old}_{\rho',1+k}\}$ for the $(1+k)$-dimensional CR model $M_{k-1}$, then according to its total nondegeneracy, one can construct the initial fields by means of the iterated Lie brackets of the generators $\mathcal L^{\sf old}_{1,1}$ and $\mathcal L^{\sf old}_{1,2}$ of $T^{1,0}M_{k-1}$ and $T^{0,1}M_{k-1}$; {\it exactly as those} for the initial vector fields on $M_k$ ({\it cf.} \thetag{\ref{iterated-bracket}})\,\,---\,\,here we assign the symbol $\sf "old"$ to objects corresponding to $M_{k-1}$. More precisely, if we have $\mathcal L_{\ell,j}=[\mathcal L_{1},\mathcal L_{\ell-1,i}]$ for $j=1,\ldots,1+k$, then correspondingly we have $\mathcal L^{\sf old}_{\ell,j}=[\mathcal L^{\sf old}_{1},\mathcal L^{\sf old}_{\ell-1,i}]$. Consequently, for a general biholomorphism $h^{\sf old}:M_{k-1}\rightarrow {\bf M}_{k-1}$ and proceeding as subsection \ref{sec-ambiguity-matrix} for the complexified push-forward $h^{\sf old}_\ast:\mathbb C\otimes TM_{k-1}\rightarrow \mathbb C\otimes T{\bf M}_{k-1}$, one finds that if ({\it cf.} Lemma \ref{lem-number-L11-12}):
\[
h_\ast({\mathcal L}_{\ell,j}):=a_1^{p}\overline a_1^q\,{\bf L}_{\ell,j}+\sum_{l<\ell}\,{\sf a}_{r_i}\,{\bf L}_{l,r} \ \ \ \ \ {\scriptstyle (j\,=\,1\,,\,\ldots\,,\,1+k),}
\]
then correspondingly we also should have:
\[
h^{\sf old}_\ast({\mathcal L}^{\sf old}_{\ell,j}):=a_1^{p}\overline a_1^q\,{\bf L}^{\sf old}_{\ell,j}+\sum_{l<\ell}\,{\sf a}_{r_i}\,{\bf L}^{\sf old}_{l,r} \ \ \ \ \ {\scriptstyle (j\,=\,1\,,\,\ldots\,,\,1+k)},
\]
though in the former case the appearing group parameter-functions are in terms of the complex variables $z',w_1',\ldots,w_{k-1}',w_k'$ and in the latter case they do not admit the last one $w_k'$. The only distinction here is that the initial frame of $M_k$ has one more initial vector field, namely $\mathcal L_{\rho,2+k}$ for which its image under $h_\ast$ should be computed, separately. This $h_\ast(\mathcal L_{\rho,2+k})$ manifests itself as the first column of $\bf g$.
\endproof

\begin{Remark}
\label{Remark-M-k-1}
By an inspection of the above proof, one finds that among the construction of the ambiguity matrix associated with $M_{k-1}$, the assigned weights to all the appearing initial vector fields, 1-forms and group parameters will be exactly as their corresponding items in the case of $M_k$.
\end{Remark}

\subsection{Associated structure equations}

According to our systematic strategy, introduced at the beginning of this section, now we are ready to compute the associated structure equations of the biholomorphic equivalence problem to the model $M_k$.
  Assuming $\Gamma:=(\Gamma_{\rho,2+k},\ldots,\Gamma_{1,1})^t$ and $\Sigma:=(\sigma_{\rho,2+k},\ldots,\sigma_{1,1})^t$ as our lifted and initial coframes, then by differentiating the the both sides of the equality \thetag{\ref{ambiguity-coframe-k}}, which can be rewritten as $\Gamma={\bf g}\cdot\Sigma$, gives:
\begin{equation}
\label{differentiation}
d\Gamma=d{\bf g}\wedge\Sigma+{\bf g}\cdot d\Sigma.
\end{equation}
For the first part $d{\bf g}\wedge\Sigma$ at the right hand side of this equation, one can replace it by:
\[
\underbrace{d{\bf g}\cdot {\bf g}^{-1}}_{\omega_{\tt MC}}\wedge \underbrace{{\bf g}\cdot\Sigma}_{\Gamma},
\]
where $\omega_{\tt MC}$ is the well-known {\sl Maurer-Cartan} matrix of the Lie group $G$. Since $\bf g$ is lower triangular with the powers of the form $a_1^r\overline a_1^s$ on its main diagonal ({\it cf.} \thetag{\ref{ambiguity-coframe-k}}), then $\omega_{\tt MC}$ is again lower triangular of the shape displaying in the following expanded form of the equation \thetag{\ref{differentiation}}:

\begin{equation}
\label{differentiation-2}
\footnotesize\aligned
\left(
  \begin{array}{c}
    d\Gamma_{\rho,i}\\
    d\Gamma_{\rho-1,j} \\
    \vdots \\
    d\Gamma_{3,5} \\
    d\Gamma_{3,4} \\
    d\Gamma_{2,3} \\
    d\Gamma_{1,2} \\
    d\Gamma_{1,1} \\
  \end{array}
\right)
&=
\underbrace{\left(
  \begin{array}{ccccc}
    p\alpha+q\overline\alpha & 0 & 0 & 0 & 0 \\
    \delta_\bullet & p'\alpha+q'\overline\alpha & 0 & 0 & 0 \\
    \vdots & \vdots & \ddots & \cdots & \cdots \\
    \delta_\bullet & \delta_\bullet & \delta_\bullet & \overline\alpha & 0 \\
    \delta_\bullet & \delta_\bullet & \delta_\bullet & \delta_\bullet & \alpha \\
  \end{array}
\right)}_{\omega_{\tt MC}}
\wedge
\left(
  \begin{array}{c}
    \Gamma_{\rho,i}\\
    \Gamma_{\rho-1,j} \\
    \vdots \\
    \Gamma_{3,5} \\
    \Gamma_{3,4} \\
    \Gamma_{2,3} \\
    \Gamma_{1,2} \\
    \Gamma_{1,1} \\
  \end{array}
\right)
\\
&
+
\underbrace{\left(
  \begin{array}{cccccccc}
    a_1^p\overline a_1^q & 0 & 0 & 0 & 0 & 0 & 0 & 0 \\
    {\sf a}_\bullet & a_1^{p'}\overline a_1^{q'} & 0 & 0 & \ldots & \ldots & 0 & 0 \\
    {\sf a}_\bullet & {\sf a}_\bullet & \ddots & 0 & 0 & \ldots & \ldots & 0 \\
    {\sf a}_\bullet & {\sf a}_\bullet & {\sf a}_\bullet & a_1\overline a_1^2 & 0 & 0 & \ldots & 0 \\
    {\sf a}_\bullet &\vdots & {\sf a}_\bullet & 0 & a_1^2\overline a_1 & 0 & \ldots & 0 \\
    {\sf a}_\bullet & \vdots & {\sf a}_\bullet & -\overline a_3 & a_3 & a_1\overline a_1 & 0 & 0 \\
    {\sf a}_\bullet & \vdots & {\sf a}_\bullet & \overline a_4 & a_5 & -\overline a_2 & \overline a_1 & 0 \\
    {\sf a}_\bullet & {\sf a}_\bullet & {\sf a}_\bullet & \overline a_5 & a_4 & a_2 & 0 & a_1 \\
  \end{array}
\right)}_{{\bf g}}
\cdot
\left(
  \begin{array}{c}
    d\sigma_{\rho,i} \\
    d\sigma_{\rho-1,j} \\
    \vdots \\
    d\sigma_{3,5} \\
    d\sigma_{3,4} \\
    d\sigma_{2,3} \\
    d\sigma_{1,2} \\
    d\sigma_{1,1} \\
  \end{array}
\right),
\endaligned
\end{equation}
with $\alpha:=\frac{d\,a_1}{a_1}$ and with $\delta_\bullet$s as some (possibly zero) certain combinations of the standard forms $da_\bullet$ with the coefficient functions in terms of $a_1,a_2,\ldots$.
The equations of \thetag{\ref{differentiation-2}} are called the {\sl structure equations} of the biholomorphic equivalence problem to $M_k$. The following lemma is encouraging enough to have some rigorous weight analysis on the structure equations in the next section. Recall that for each term $a_jd\sigma_{\ell,i}$, coming from the last matrix multiplication of \thetag{\ref{differentiation-2}}, the associated weight is naturally defined  as $[a_j]+[d\sigma_{\ell,i}]$.

\begin{Lemma}
\label{lemma-weight-adsigma}
All entries of the last vertical matrix ${\bf g}\cdot d\Sigma$ at the right hand side of the above structure equations \thetag{\ref{differentiation-2}} are homogeneous of the equal weight zero.
\end{Lemma}

\proof
It is a straightforward consequence of Lemma \ref{lemma-weight}, reminding that the assigned weight to each $\sigma_\ell$ and its differentiation $d\sigma_\ell$ is $-\ell$.
\endproof

\subsection{Torsion coefficients}
\label{subsec-torsion-coef}

Our next aim is to restate the above structure equations \thetag{\ref{differentiation-2}} absolutely independent of the initial 1-forms $\sigma_{\ell,i}$ and their differentiations. For this purpose, we shall focus on the second matrix term ${\bf g}\cdot d\Sigma$. The Darboux-Cartan structure computed in Proposition \ref{Darboux-Cartan} enables one to replace each 2-form $d\sigma_\bullet$ by some combination of the wedge products between initial 1-forms $\sigma_\bullet$. Afterward, by means of the equality $\Sigma={\bf g}^{-1}\cdot\Gamma$, it is also possible to replace each initial 1-form $\sigma_\bullet$ by some combination of the lifted 1-forms $\Gamma_\bullet$. Doing so, then all differentiations at the right hand side vertical matrix ${\bf g}\cdot d\Sigma$ of \thetag{\ref{differentiation-2}} will be expressible in terms of the wedge products of the lifted 1-forms $\Gamma_\bullet$. Consequently, our structure equations will be converted into the form:
\begin{equation}
\label{structure-equations-old}
\aligned
d\,\Gamma_{\ell,i}&:=(p_i\,\alpha+q_i\,\overline\alpha)\wedge\Gamma_{\ell,i}+\sum_{r,j,\  l\gneqq\ell}\delta_r\wedge\Gamma_{l,j}
\\
&+\sum_{l,j,m,n}\,T^{i}_{jn}(a_\bullet)\,\Gamma_{l,j}\wedge\Gamma_{m,n}, \ \ \ \ \ \ \ \ \ \ \ \ \ \ \ \ \ \ \ \ \ {\scriptstyle (\ell\,=\,1\,,\,\ldots\,,\,\rho, \ \ i\,=\,1\,,\,\ldots\,,\,2+k)},
\endaligned
\end{equation}
where $T^{i}_{jn}$s are some certain functions in terms of the group parameters $a_\bullet$ which are called by the {\sl torsion coefficients} of the problem.

\begin{Remark}
\label{remark-g-1}
Since our ambiguity matrix $\bf g$ is invertible and lower triangular with the powers $a_1^p\overline a_1^q$ at its diagonal, then a simple induction on the number of its column and rows shows that ${\bf g}^{-1}$ is again lower triangular where its non-diagonal entries are some fraction polynomial functions with some powers of the form $a_1^r\overline a_1^s$ as their denominators. Also, if the $i$-th diagonal entry of $\bf g$ is, say, $a_1^p\overline a_1^q$ then this entry in ${\bf g}^{-1}$ is $\frac{1}{a_1^p\overline a_1^q}$. Finally, thanks to Lemma \ref{lem-number-L11-12} and again since $\bf g$ is lower triangular, then in the expression of each length $\ell$ lifted 1-from $\Gamma_{\ell,i}$ as \thetag{\ref{ambiguity-coframe-k}}, the only appearing initial 1-form of the lengths $\leqslant\ell$ is $\sigma_{\ell,i}$. Consequently, by a backward induction on the length $\ell$ of the initial 1-forms from $\rho$ to $1$, we discover a same fact in expressing each initial 1-form $\sigma_{\ell,i}$ in terms of the lifted ones through the equality $\Sigma={\bf g}^{-1}\cdot\Gamma$: the only appearing lifted 1-form in the expression of $\sigma_{\ell,i}$ of the length $\leqslant\ell$ is $\Gamma_{\ell,i}$. Therefore, if the $j$-th row of the vertical matrix $\Sigma$, say  $\sigma_{\ell,i}$, is of the length $\ell$ then the $j$-th row of ${\bf g}^{-1}$ is of the form:
\[
\big(c_\bullet,\ldots,c_\bullet, \underbrace{0, \ldots, 0}_{t_1 \textrm{times}}, \underbrace{\frac{1}{a_1^r\overline a_1^s}}_{i-th \ \textrm{ place }}, \underbrace{0,\ldots,0}_{t_2 \ \textrm{times}}\big)
\]
where $r+s=\ell$ and $t_1+t_2+1$ is equal or more than the number of initial 1-forms $\sigma_\bullet$ of the lengths $\leqslant\ell$.
\end{Remark}

\section{Weight analysis on the structure equations}
\label{section-weight-analysis}

In the previous section, we assigned naturally some weights to the complex variables, initial and lifted vector fields and 1-forms, their differentiations and also to group parameters. The main purpose of this section is to show that all the appearing torsion coefficients in the constructed structure equations \thetag{\ref{structure-equations-old}} are weighted homogeneous of the same weight zero.  For this aim, we inspect more the structure of the inverse matrix ${\bf g}^{-1}$ via some auxiliary lemmas. But at first we need the following definition.

\begin{Definition}
Let:
\[
f(a_1, a_2,\ldots)=\frac{a_1^{r_1}\overline a_1^{s_1}a_2^{s_2}\overline a_2^{s_2}\ldots a_n^{r_n}\overline a_n^{s_n}}{a_1^{r}\overline a_1^s}
\]
be an arbitrary monomial fraction in terms of the group parameters. Then, the weight of $f$ is defined as:
\[
[f]=r_1[a_1]+s_1[\overline a_1]+r_2[a_2]+s_2[\overline a_2]+\ldots+r_n[a_n]+s_n[\overline a_n]-r[a_1]-s[\overline a_1].
\]
A weighted homogeneous polynomial fraction  is a sum of monomial fractions of the same weigh.
\end{Definition}

As stated in Lemma \ref{lemma-weight}, all the nonzero entries in a fixed column of our ambiguity matrix $\bf g$ are of the same weight. Our next goal is to show that in the inverse matrix ${\bf g}^{-1}$, the {\it rows} enjoy a similar fact.

\begin{Lemma}
\label{lemma-weight-g-1}
Fix an integer $i_0=1,\ldots,2+k$ and let $-\ell$ be the weight of a certain $1$-form $\sigma_\ell$ standing at the $i_0$-th row of the vertical matrix $\Sigma=\big(\sigma_{\rho,2+k}, \ldots, \sigma_{1,2}, \sigma_{1,1}\big)^t$
in \thetag{\ref{ambiguity-coframe-k}}. Then,
\begin{itemize}
\item[(i)] all the nonzero entries of the $i_0$-th row of ${\bf g}^{-1}$ are of the same homogeneous weight $-\ell$, too.
\item[(ii)] if the $j$-th row of $\Sigma$ is of the weight $-(\ell+1)$ and if the $(i_0j)$-th entry of $\bf g$ is  $e_{i_0j}$ then, this entry in ${\bf g}^{-1}$ is of the form:
\[
-\frac{e_{i_0j}}{a_1^m\overline a_1^n},
\]
for some constant integers $m$ and $n$.
\end{itemize}
\end{Lemma}

\proof
We prove the both parts by an induction on the codimension $k$ of the models. The base of this induction is provided by inspecting the matrices introduced in \cite[p. 89]{Samuel} for $k=2$ and \cite[p. 104]{5-cubic} for $k=3$\,\,---\,\,according to the Conjecture \ref{conjecture} we are considering CR models of the lengths $\rho\geqslant 3$ which start from $k=2$.  Assume that the assertions hold for all CR models of codimensions $<k-1$.  By Proposition \ref{prop-g-k-1}, if ${\bf g}_{k-1}$ is the ambiguity matrix of the equivalence problem to the CR model $M_{k-1}$, then:

 \begin{equation}
 \label{g-inverse-k-1}
 \footnotesize\aligned
 {\bf g}^{-1}=\left(
  \begin{array}{cc}
    \frac{1}{a_1^{p}\overline a_1^{q}} & 0 \,\,\,\, 0 \,\,\,\, \cdots \,\,\,\, 0 \,\,\,\, 0  \\
    \begin{array}{c}
      {\bf 0} \\
      {\sf b}_j \\
      \vdots \\
     {\sf b}_1
    \end{array}
         &
     \left(
       \begin{array}{ccccc}
          &  &  &  &    \\
          &  &  &  &    \\
          &  & {\bf g}^{-1}_{k-1}  &  &  \\
          &  &  &  &    \\
       \end{array}
     \right)
           \\
  \end{array}
\right),
 \endaligned
 \end{equation}
 with $p+q=\rho$ and
 for some certain functions ${\sf b}_\bullet$. Thus, according to our induction and by Remark \ref{Remark-M-k-1}, it suffices to prove (i) just for each entry ${\sf b}_t$ at some $i_0$-th row of ${\bf g}^{-1}$ with $i_0\neq 1$. According to Lemma \ref{lemma-weight}, all the nonzero group parameters at the first column of $\bf g$ are of the same maximum weight $\rho$. By our induction hypothesis and except ${\sf b}_t$, we know that all the nonzero entries at the $i_0$-th row of ${\bf g}^{-1}$ are of the same weight $-\ell$. We show that if ${\sf b}_t\neq 0$, then it is of the same weight, too. Multiplying the $i_0$-th row of ${\bf g}^{-1}$ by the first column of $\bf g$ gives:
 \[
 {\sf b}_t\cdot (a_1^{p}\overline a_1^{q})+\Psi=0
 \]
where $\Psi$ is some function of the weight $\rho-\ell$. Taking into account that $p+q=\rho$, then the polynomial fraction ${\sf b}_t=-\frac{\Psi}{a_1^{p}\overline a_1^{q}}$ is of the homogeneous weight $\rho-\ell-\rho=-\ell$, as was expected.

\noindent
For the second part (ii), and according to our induction, it suffices to prove it only for some of the entries $e_{i_01}$ at the first column of $\bf g$, namely for $j=1$. Since the first row of $\Sigma$ is of the weight $-\rho$, then we have to look for weight $-\ell=-(\rho-1)$ rows $i_0$ of the inverse matrix ${\bf g}^{-1}$. By the first part (i), these rows are in front of the weight $-(\rho-1)$ initial 1-forms $\sigma_{\rho-1,i}$ in the equation $\Sigma={\bf g}^{-1}\cdot\Gamma$ and hence $e_{i_01}$ stands below the zero vector $\bf 0$ at the first column of  $\bf g$ ({\it cf.} \thetag{\ref{ambiguity-coframe-k}}). Hence, $i_0\neq 1$. Assume that ${\sf b}_r$ stands at the same entry of ${\bf g}^{-1}$ as $e_{i_01}$ in $\bf g$. We aim to show ${\sf b}_r=-\frac{e_{i_01}}{a_1^m\overline a_1^n}$. The $i_0$-th row of ${\bf g}^{-1}$ is of the form ({\it cf.} Remark \ref{remark-g-1}):
\[
\big({\sf b}_r,c_1,\ldots,c_t,0,\ldots,0,\underbrace{\frac{1}{a_1^r\overline a_1^s}}_{\textrm{ $i_0$-th place}}, 0,\ldots,0\big),
\]
where $t+1$ is the number of the weight $-\rho$ lifted 1-forms $\sigma_\rho$. Then, multiplying again the above $i_0$-th row of ${\bf g}^{-1}$ to the first column of $\bf g$ and granted the Remark \ref{Remark-M-k-1} about the zero vector $\bf 0$ at this column gives:
\[
\big({\sf b}_r,c_1,\ldots,c_t,0,\ldots,0,\underbrace{\frac{1}{a_1^r\overline a_1^s}}_{\textrm{ $i_0$-th place}}, 0,\ldots,0\big)\cdot\big(a_1^p\overline a_1^q, \underbrace{{\bf 0}}_{\textrm{ $t$-tuple}}, \ldots, \underbrace{e_{i_01}}_{\textrm{ $i_0$-th place}}, \ldots\big)^t=0.
\]
Now, simplifying this equality after multiplication and solving it in terms of ${\sf b}_r$ gives ${\sf b}_r=-\frac{e_{i_01}}{a_1^{p+r}\overline a_1^{q+s}}$,
as desired.
\endproof

Roughly speaking, the first part $(i)$ of this lemma states that for each fixed row of the three matrices appearing in the equation $\Sigma={\bf g}^{-1}\cdot\Gamma$, all the nonzero entries are of the same negative weight. Furthermore, taking into account the shape of the lower triangular matrix ${\bf g}^{-1}$ and by the first part of the above lemma, one observes that ({\it see} also Remark \ref{remark-g-1});

\begin{Lemma}
\label{sigma-Gamma-expression}
For each weight $-\ell$ initial $1$-form $\sigma_{\ell,i}$, its expression in terms of the lifted $1$-forms is as follows:
\[
\sigma_{\ell,i}:=\sum_{l\gneqq\ell}\,{\sf A}^i_j(a_\bullet)\,\Gamma_{l,j}+\frac{1}{a_1^{p_i}\overline a_1^{q_i}}\,\Gamma_{\ell,i},
\]
with $p_i+q_i=\ell$ and for some weighted homogeneous polynomial fractions ${\sf A}^i_j$ of the weight $-\ell$ where their denominators are some powers of only $a_1$ and $\overline a_1$.
\end{Lemma}

Now, we are ready to prove the main result of this section;

\begin{Proposition}
\label{prop-torsion-weight-zero}
All torsion coefficients $T^i_{jn}(a_\bullet)$ appearing among the structure equations \thetag{\ref{structure-equations-old}} are weighted homogeneous polynomial fractions of the equal weight zero where their denominators are some powers of only $a_1$ and $\overline a_1$.
\end{Proposition}

\proof
According to \thetag{\ref{differentiation-2}}, each structure equation can be expressed as:
\[
d\Gamma_{\ell,i}=(p_i\alpha+q_i\overline\alpha)\wedge\Gamma_{\ell,i}+\sum_{l\gneqq\ell}\,\delta_{i_j}\wedge\Gamma_{l,j}+\sum_{l\gneqq\ell}\,a_{i_j}d\sigma_{l,j}+a_1^{p_i}\overline a_1^{q_i}\,d\sigma_{\ell,i}
\]
with $p_i+q_i=\ell$.
Our torsion coefficients come from the last parts:
\begin{equation}
\label{torsion}
\sum_{l\gneqq\ell}\,a_{i_j}d\sigma_{l,j}+a_1^{p_i}\overline a_1^{q_i}\,d\sigma_{\ell,i}
\end{equation}
of this equation after replacing each differentiation $d\sigma_\bullet$ according to the Darboux-Cartan structure computed in Proposition \ref{Darboux-Cartan} and next substituting each initial 1-form $\sigma_\bullet$ with some combinations of lifted 1-forms $\Gamma_\bullet$ by means of the equality $\Sigma={\bf g}^{-1}\cdot\Gamma$. Thanks to Lemma \ref{lemma-weight-adsigma}, the weight of the coefficient $a_{i_j}$ in the term $a_{i_j}d\sigma_{l,j}$ of \thetag{\ref{torsion}} is $l$ . Moreover, according to Proposition \ref{Darboux-Cartan} we have:
\[
d\sigma_{l,j}:=\sum_{\beta,\gamma}\,{\sf c}_{\beta,\gamma}\,\sigma_\beta\wedge\sigma_\gamma \ \ \ \ {\rm with} \ \ \ \ \beta+\gamma=l.
\]
After replacing the expressions of $\sigma_\beta$ and $\sigma_\gamma$ as Lemma \ref{sigma-Gamma-expression}, such Darboux-Cartan structure takes the form:
\[
a_{i_j}\,d\sigma_{l,j}=\sum_{l_1+l_2\geqslant l}\,\big (a_{i_j}\,{\sf T}^j_{m,n}(a_\bullet)\big )\,\Gamma_{l_1,m}\wedge\Gamma_{l_2,n}
\]
where the polynomial fractions ${\sf T}^j_{m,n}$ are multiplications of some weight $-\beta$ and $-\gamma$ polynomial fractions with $\beta+\gamma=l$. Thus, all the coefficients ${\sf T}^j_{m,n}$ are of the same weight $-\beta-\gamma=-l$ and hence, each coefficient $a_{i_j}{\sf T}^j_{m,n}(a_\bullet)$ in the above expression is of the weight zero. Similar fact holds true also for the last term $a_1^{p_i}\overline a_1^{q_i}\,d\sigma_{\ell,i}$ of \thetag{\ref{torsion}}. Now, each torsion coefficient $T^i_{m,n}$ is made as the sum of coefficients of $\Gamma_{l_1,m}\wedge\Gamma_{l_2,n}$ in the expressions of all terms $a_{i_j}d\sigma_{l,j}$ and $a_1^{p_i}\overline a_1^{q_i}\,d\sigma_{\ell,i}$, visible in \thetag{\ref{torsion}}. Therefore, it is of the weight zero, as claimed. The second part of the assertion is a consequence of Remark \ref{remark-g-1}.
\endproof

Before concluding this section, let us present another result of the second part (ii) of Lemma \ref{lemma-weight-g-1};

\begin{Lemma}
\label{l-l-1}
 If in the structure equation $d\Gamma_{\ell-1,m}$ of \thetag{\ref{differentiation-2}} we have the term $a_jd\sigma_{\ell,n}$ for some (possibly zero) group parameter $a_j$, then the coefficient of $\Gamma_{\ell,n}$ in the expression of $\sigma_{\ell-1,m}$, through the equation $\Sigma={\bf g}^{-1}\cdot\Gamma$, is of the form $-\frac{a_j}{a_1^{r}\overline a_1^s}$ for some constant integers $r$ and $s$.
\end{Lemma}

\proof
 The term $a_jd\sigma_{\ell,n}$ in \thetag{\ref{differentiation-2}} comes only from the second part ${\bf g}\cdot d\Sigma$ and hence the appearance of this term in the structure equation $d\Gamma_{\ell-1,m}$ means that the coefficient of $\sigma_{\ell,n}$ in the expression of $\Gamma_{\ell-1,m}$\,\,---\,\,coming from the equality $\Gamma={\bf g}\cdot\Sigma$\,\,---\,\,is $a_j$:
\[
\footnotesize\underbrace{\left(
  \begin{array}{c}
    \vdots \\
    \vdots \\
    \Gamma_{\ell-1,m} \\
    \vdots \\
  \end{array}
\right)}_{\Gamma}
=
\underbrace{\left(
  \begin{array}{cccc}
      \vdots & \vdots & \vdots & \vdots \\
      \vdots & \vdots & \vdots & \vdots \\
    \ldots & a_j & \ldots & \ldots \\
    \vdots & \vdots & \vdots & \vdots \\
  \end{array}
\right)}_{\bf g}
\cdot
\underbrace{\left(
  \begin{array}{c}
    \vdots \\
    \sigma_{\ell,n} \\
    \vdots \\
    \vdots \\
  \end{array}
\right).}_{\Sigma}
\]
By the above matrix equation and according to the second part (ii) of Lemma \ref{lemma-weight-g-1}, we will have some $-\frac{a_j}{a_1^{r}\overline a_1^s}$ in ${\bf g}^{-1}$ in place of the same entry $a_j$ in $\bf g$. But this entry in the inverse matrix determines, through the equality $\Sigma={\bf g}^{-1}\cdot\Gamma$, the coefficient of the lifted 1-form $\Gamma_{\ell,n}$ in the expression of $\sigma_{\ell-1,m}$.
\endproof

This suggests that if we are seeking the coefficient of $\Gamma_{\ell,n}$ in the expression of some $\sigma_{\ell-1,m}$, then it is opposite to the fraction of the coefficient of $d\sigma_{\ell,n}$ in the structure equation $d\Gamma_{\ell-1,m}$ by some powers of $a_1$ and $\overline a_1$. This result will be of much use in the next section.

\section{Picking up an appropriate weighted homogeneous subsystem}
\label{section-picking-system}

Now we are ready to apply Cartan's method on the biholomorphic equivalence problem of the CR model $M_k$. The first two essential steps of this method are {\it absorbtion} and {\it normalization}, based on some fundamental results introduced in \cite[Proposition 4.7]{5-cubic} ({\it see also \cite{Olver-1995}}). According to these results, one is permitted to substitute as follows each Maurer-Cartan 1-form $\alpha$ and $\delta_j$ in the structure equations \thetag{\ref{structure-equations-old}}:
\begin{equation}
\label{replacement-maurer-cartan}
\aligned
\alpha&\mapsto \alpha+t_{2+k}\,\Gamma_{\rho,2+k}+\ldots+t_2\, \Gamma_{1,2}+t_1\,\Gamma_{1,1},
\\
\delta_j&\mapsto \delta_j+s^j_{2+k}\,\Gamma_{\rho,2+k}+\ldots+s^j_2\, \Gamma_{1,2}+s^j_1\,\Gamma_{1,1},
\endaligned
\end{equation}
for {\it arbitrary} coefficient functions $t_\bullet$ and $s^\bullet_\bullet$. We can apply such substitutions and try to convert new (torsion) coefficients of the wedge products $\Gamma_{\ell_1,i_1}\wedge\Gamma_{\ell_2,i_2}$ to some constant integers\,\,---\,\,possibly zero\,\,---\,\, by appropriate determinations of the {arbitrary} functions $t_\bullet, s^\bullet_\bullet$ (this is the absorption step).
For this purpose, it may be inadequate only such determination of these arbitrary functions but it necessitates also to determine\,\,---\,\,or normalize in this literature\,\,---\,\,some of the group parameters, appropriately in terms of the other ones by equating to zero (or other constants) still remaining non-constant coefficients. These coefficients are called the {\it essential torsion coefficients}.

 Thus to proceed along the absorption and normalization steps, one has to solve an arising polynomial system with $t_\bullet$, $s^\bullet_\bullet$ and some of the group parameters as its unknowns. The virtual importance of the solution of this system is not determining the coefficient functions $t_\bullet$ and $s^\bullet_\bullet$ but it is actually the found values of involving group parameters $a_\bullet$. Unfortunately, solving such arising polynomial system, specifically in this general manner, causes certainly some unavoidable and serious algebraic complexity. The main purpose of this section is to bypass and manipulate such complexity by picking up an appropriate and convenient subsystem that affords to bring all results we are seeking from the solution of the original system. Before explaining our practical method of constructing this desired susbsystem\,\,---\,\,which will be divided into two major parts\,\,---\,\,at first, we need the following auxiliary lemma;

 \begin{Lemma}
 \label{lem-struc-unique}
   Assume that $\sigma_{\ell-1,i}\wedge\sigma_{1,t}$, for $t=1$ or $2$, is the unique appearing wedge product in the Darboux-Cartan structure of $d\sigma_{\ell,j}$, as stated in Lemma \ref{lem-unique-1-form}. Then, among all the expressions of differentiations $d\sigma_{l,r}$, with $l\geqslant\ell$, in terms of the wedge products of the lifted $1$-forms, a nonzero coefficient of $\Gamma_{\ell-1,i}\wedge\Gamma_{1,t}$ appears uniquely in $d\sigma_{\ell,j}$. Such coefficient is a fraction of the form $\frac{1}{a_1^p\overline a_1^q}$ for some constant integers $p$ and $q$.
 \end{Lemma}
 \proof
By Remark \ref{remark-g-1}, in the expression of each $\sigma_{\ell',r}$ through the equality $\Sigma={\bf g}^{-1}\cdot\Gamma$, the only appearing lifted 1-form $\Gamma_{l,m}$ with $l\leqslant\ell'$ is some $\frac{1}{a_1^p\overline a_1^q}\Gamma_{\ell',r}$. In particular, the only initial 1-form having some coefficient of $\Gamma_{1,1}$ in its expression is $\sigma_{1,1}$ and this coefficient is $\frac{1}{a_1}$. Similarly, the only initial 1-form having some coefficient of $\Gamma_{1,2}$ is $\sigma_{1,2}$ with the coefficient $\frac{1}{\overline a_1}$. Consequently, in the expression of a fixed differentiation $d\sigma_{l_0,r}$ with $l_0\geqslant\ell$, one finds some nonzero coefficient of $\Gamma_{\ell-1,i}\wedge\Gamma_{1,t}$ whenever in its Darboux-Cartan structure, $d\sigma_{l_0,r}$ includes some nonzero coefficient of the wedge product $\sigma_{l',j}\wedge\sigma_{1,t}$ with $l'\leqslant\ell-1$. We claim that $\sigma_{l',j}=\sigma_{\ell-1,i}$. Since $\sigma_{l',j}\wedge\sigma_{1,t}$ appears in the Darboux-Cartan structure of $d\sigma_{l_0,r}$ then Proposition \ref{Darboux-Cartan} implies that $l'+1=l_0\geqslant\ell$ and whence $l'\geqslant\ell-1$. Consequently, $l'=\ell-1$ and thus  $\sigma_{l',j}=\sigma_{\ell-1,j}$. Furthermore, again by Remark \ref{remark-g-1}, $\sigma_{\ell-1,i}$ is the only weight $-(\ell-1)$ initial form containing some nonzero coefficient of $\Gamma_{\ell-1,i}$ in its expression. This results that $\sigma_{\ell-1,j}=\sigma_{\ell-1,i}$, as was claimed. But on the other hand, according to our assumption, $\sigma_{\ell-1,i}\wedge\sigma_{1,t}$ appears uniquely in the Darboux-Cartan structure of $d\sigma_{\ell,j}$ and hence we should have $d\sigma_{l_0,r}=d\sigma_{\ell,j}$, as was desired. In addition, the coefficient $\Gamma_{\ell-1,i}\wedge\Gamma_{1,t}$ in $d\sigma_{\ell,j}$ comes from the wedge product $\sigma_{\ell-1,i}\wedge\sigma_{1,t}$ in its Darboux-Cartan structure and by what mentioned at the beginning of the proof, it will be nothing but some fraction $\frac{1}{a_1^p\overline a_1^q}$.
\endproof

\subsection{Picking up an appropriate subsystem}
\label{sub-pick}

Our strategy of picking up appropriate torsion coefficients from the structure equations, after absorption, is divided into two essential parts depending upon the weight.

\subsubsection{First part: structure equations of the weights $-\ell\,=\,-1, \ldots, -(\rho-1)$}

Consider:
 \begin{equation}
 \label{struc-eq-ell}
\aligned
d\Gamma_{\ell,m}=(p_m\alpha+q_m\overline\alpha)\wedge\Gamma_{\ell,m}+\sum_{l\gneqq\ell}\,\delta_{i_t}\wedge\Gamma_{l,j}
+\sum_{l\geq\ell +2 }\,a_{j_n}d\sigma_{l,n}+\sum_{r}\,a_{j_r}d\sigma_{\ell+1,r}+a_1^{p_m}\overline a_1^{q_m}\,d\sigma_{\ell,m},
\endaligned
\end{equation}
 as a weight $-\ell$ structure equation in \thetag{\ref{differentiation-2}} for $\ell=1, \ldots, \rho-1$. {\it We focus just on the terms} $a_{j_{r}}d\sigma_{\ell+1,r}$ in the penultimate part $\sum_{r}\,a_{j_r}d\sigma_{\ell+1,r}$ of this structure equation. Lemma \ref{lemma-weight-adsigma} implies that the group parameters $a_{j_r}$, visible in it, are of the weight $\ell+1$. As a consequence of the above Lemma \ref{lem-struc-unique} and in the expression of each {\it fixed} term $a_{j_{r_0}}d\sigma_{\ell+1,r_0}$, in terms of the wedge products of lifted 1-forms, one finds a certain product:
 \begin{equation}
 \label{1}
 \frac{a_{j_{r_0}}}{a_1^{p_\bullet}\overline a_1^{q_\bullet}}\,\Gamma_{\ell,i_j}\wedge\Gamma_{1,t_r} \ \ \ \ {\scriptstyle (t_r\,=1  \ {\rm or} \ 2)},
 \end{equation}
 coming from some $\sigma_{\ell,i_j}\wedge\sigma_{1,t_r}$,  uniquely appearing in the Darboux-Cartan structure of $d\sigma_{\ell+1,r_0}$, such that no any other term in the part $\sum_{l\geq\ell +2, n}\,a_{j_n}d\sigma_{l,n}+\sum_{r\neq r_0}\,a_{j_r}d\sigma_{\ell+1,r}$ of \thetag{\ref{struc-eq-ell}} brings any nonzero coefficient of  it. Then, as is our strategy, we seek for all coefficients of this wedge product $\Gamma_{\ell,i_j}\wedge\Gamma_{1,t_r}$ in \thetag{\ref{struc-eq-ell}}. Let us do it part by part.

We continue with the last term $a_1^{p_{m}}\overline a_1^{q_{m}}\,d\sigma_{\ell,m}$. Assuming the Darboux-Cartan structure:
\[
 d\sigma_{\ell,m}=\sum_{\ell_1+\ell_2=\ell}{\sf c}_{\ell_1,\ell_2}\sigma_{\ell_1,t}\wedge\sigma_{\ell_2,s},
\]
then the desired wedge product $\Gamma_{\ell,i_j}\wedge\Gamma_{1,t_r}$ is producible only by the terms of the form\footnote{Remind that one can find the lifted 1-forms $\Gamma_{1,1}$ and $\Gamma_{1,2}$ only in the expressions of $\sigma_{1,1}$ and $\sigma_{1,2}$, respectively.} ${\sf c}_{t}\,\sigma_{\ell-1,t}\wedge\sigma_{1,t_r}$. In order to find the coefficient of this product, we have to pick the coefficient of $\Gamma_{\ell,i_j}$ in the expression of $\sigma_{\ell-1,t}$s. By Lemma \ref{l-l-1} and if the coefficient of $d\sigma_{\ell,i_j}$ in the structure equation of $d\Gamma_{\ell-1,t}$ is a weight $\ell$ group parameter $a_{i_t}$, then the desired coefficient of $\Gamma_{\ell,i_j}$ in $\sigma_{\ell-1,t}$ will be of the form $-\frac{a_{i_t}}{a_1^{r'_\bullet}\overline a_1^{s'_\bullet}}$ for some constant integers $r'_\bullet$ and $s'_\bullet$\footnote{Notice also that here $a_{i_t}$ can be zero and it does not effect our next results.}. Hence the last term $a_1^{p_{m}}\overline a_1^{q_{m}}\,d\sigma_{\ell,m}$ may produce some term like:
  \begin{equation}
  \label{2}
  -\big(\sum_t\,{\sf c_t}\frac{a_{i_t}}{a_1^{r_\bullet}\overline a_1^{s_\bullet}}\big)\,\Gamma_{\ell,i_j}\wedge\Gamma_{1,t_r},
  \end{equation}
  after simplification. For later use, we emphasize from the above procedure that;

  \begin{Lemma}
  \label{lem-a-i-t}
  For each weight $\ell$ group parameter $a_{i_t}$ appearing in \thetag{\ref{2}}, there is a term $a_{i_t}d\sigma_{\ell,i_j}$ in some weight $-(\ell-1)$ structure equations $d\Gamma_{\ell-1,t}$.
  \end{Lemma}

Obviously, the second term $\sum_{l\gneqq\ell}\,\delta_{i_t}\wedge\Gamma_{l,j}$ of \thetag{\ref{struc-eq-ell}} will not produce any nonzero coefficient of $\Gamma_{\ell,i_j}\wedge\Gamma_{1,t_r}$ while after the substitutions \thetag{\ref{replacement-maurer-cartan}} in the first part $(p_m\alpha+q_m\overline\alpha)\wedge\Gamma_{\ell,m}$, one may find some terms like:
\begin{equation}
 \label{3'}
 -(p_{m}t_1+q_{m}\overline t_2)\,\Gamma_{\ell,m}\wedge\Gamma_{1,1} \ \ \ \ \textrm{ and} \ \ \ \ -(p_{m}t_2+q_{m}\overline t_1)\,\Gamma_{\ell,m}\wedge\Gamma_{1,2},
 \end{equation}
where in the case that $m=i_j$, one of them will be the sought product $\Gamma_{\ell,i_j}\wedge\Gamma_{1,t_r}$, in question.

Then, all possible coefficients of the wedge product $\Gamma_{\ell,i_j}\wedge\Gamma_{1,t_r}$ in the above weight $-\ell$ structure equation \thetag{\ref{struc-eq-ell}}, after absorption, are those presented in \thetag{\ref{1}}--\thetag{\ref{3'}}. Equating this coefficient to zero\,\,---\,\,as is the method of absorption-normalization\,\,---\,\,then one finds some fraction polynomial equation of the form:
\begin{equation*}
 \frac{a_{j_{r_0}}}{a_1^{p_\bullet}\overline a_1^{q_\bullet}}-\sum_t\,{\sf c_t}\frac{a_{i_t}}{a_1^{r_\bullet}\overline a_1^{s_\bullet}}={\sf a}_{j_{r_0}}t_1+{\sf b}_{j_{r_0}}t_2+{\sf a}_{j_{r_0}}'\overline t_1+{\sf b}_{j_{r_0}}'\overline t_2,
 \end{equation*}
for some (possibly zero) constants ${\sf a}_{j_{r_0}}, {\sf b}_{j_{r_0}}, {\sf a}_{j_{r_0}}', {\sf b}_{j_{r_0}}'$. The left hand side of this equation is actually the torsion coefficient $T^m_{i_j,t_r}$ of $\Gamma_{\ell,i_j}\wedge\Gamma_{1,t_r}$ in the structure equation \thetag{\ref{struc-eq-ell}} which comes from \thetag{\ref{1}} and \thetag{\ref{2}}. Hence according to Proposition \ref{prop-torsion-weight-zero}, it is of the weight zero. Minding that here $a_{j_{r_0}}$ is a weight $\ell+1$ group parameter while $a_{i_t}$s are of the weight $\ell$, then multiplying both side of this equation by the denominator $a_1^{p_\bullet}\overline a_1^{q_\bullet}$ gives the following equivalent {\it weighted homogeneous} polynomial equation\,\,---\,\,here we assign naturally the weight zero to the parameters $t_1$, $t_2$ and their conjugations:
\begin{equation}
\label{syst-ell}
\boxed{
a_{j_{r_0}}-\sum_t\,{\sf c}_t\,a_1^{r'_\bullet}\overline a_1^{s'_\bullet} a_{i_t}=a_1^{p_\bullet}\overline a_1^{q_\bullet}\big({\sf a}_{j_{r_0}}t_1+{\sf b}_{j_{r_0}}t_2+{\sf a'}_{j_{r_0}}\overline t_1+{\sf b'}_{j_{r_0}}\overline t_2\big)}.
\end{equation}

 \begin{Proposition}
Let $\ell=1,\ldots,\rho-1$. Then, among the procedure of absorbtion and associated with each weight $\ell+1$ group parameter $a_{j_{r_0}}$ appearing in an arbitrary weight $-\ell$ structure equation \thetag{\ref{struc-eq-ell}}, one finds a weighted homogeneous parametric complex polynomial equation as \thetag{\ref{syst-ell}}, expressing $a_{j_{r_0}}$ in terms of $a_1, \overline a_1$, some weight $\ell$ group parameters $a_{i_t}$, two parameters $t_1, t_2$ and their conjugations.
 \end{Proposition}

Let us denote by $\sf S$ the weighted homogeneous system of equations mentioned in the above proposition. Notice that $\sf S$ does not involve necessarily all the group parameters $a_\bullet$. Importantly in this system, if there is an equation like \thetag{\ref{syst-ell}} that expresses a weight $\ell+1$ group parameter $a_{j_{r_0}}$ in terms of the weight $\ell$ parameters $a_{i_t}$, then Lemma \ref{lem-a-i-t} guarantees that also for each $a_{i_t}$, we find another equation expressing it in terms of the lower weight group parameters when we perform the above method for weight $-(\ell-1)$ structure equations.

\begin{Proposition}
\label{prop-S}
For each group parameter $a_j\neq a_1$ visible in $\sf S$, there exists some weighted homogeneous equation expressing it in terms of some lower weight group parameters $a_\bullet$ and the parameters $t_1, t_2$.
\end{Proposition}

   Among the system $\sf S$, two equations coming from the first two structure equations:
 \[
 \aligned
 d\Gamma_{2,3}&=(\alpha+\overline\alpha)\wedge\Gamma_{2,3}+\sum_{l\gneqq 2}\,\delta_{i_j}\wedge\Gamma_{l,j}+\sum_{l\gneqq 3}\,{a_{i_j}}\,d\sigma_{l,j}+a_3\,d\sigma_{3,4}+\overline a_3\,d\sigma_{3,5}+a_1\overline a_1\, d\sigma_{2,3},
 \\
  d\Gamma_{1,1}&=\alpha\wedge\Gamma_{1,1}+\sum_{l\gneqq 1}\,\delta_{i_j}\wedge\Gamma_{l,j}+\sum_{l\gneqq 2}\,{a_{i_j}}\,d\sigma_{l,j}+a_2\,d\sigma_{2,3}+a_1\, d\sigma_{1,1}
 \endaligned
 \]
 are of particular importance. According to our suggested method, in the weight $-2$ structure equation $d\Gamma_{2,3}$ we should focus on the term $a_3\,d\sigma_{3,4}$ since $d\sigma_{3,4}$, together with $d\sigma_{3,5}$, are the only weight $-(2+1)=-3$ differentiations visible in it. Since $\mathcal L_{3,4}=[\mathcal L_{1,1}, \mathcal L_{2,3}]$, then the uniquely appearing wedge product in the Darboux-Cartan structure of $d\sigma_{3,4}$ is $\sigma_{2,3}\wedge\sigma_{1,1}$ ({\it cf.} Lemma \ref{lem-unique-1-form} and its proof). Thus, we shall look for the (torsion) coefficient of $\Gamma_{2,3}\wedge\Gamma_{1,1}$ in this structure equation $d\Gamma_{2,3}$. Also in the weight $-1$ structure equation $d\Gamma_{1,1}$ we should focus on the single term $a_2\,d\sigma_{2,3}$. The uniquely appearing wedge product in the Darboux-Cartan structure of $d\sigma_{2,3}$ is $\sigma_{1,2}\wedge\sigma_{1,1}$, then let us find the coefficient of $\Gamma_{1,2}\wedge\Gamma_{1,1}$ in this structure equation. Performing necessary computations, we respectively find the following two weight zero homogeneous equations, after applying the substitutions \thetag{\ref{replacement-maurer-cartan}}:
\begin{equation}
\label{t1-2}
\aligned
\frac{a_3}{a_1^2\overline a_1}+i\frac{\overline a_2}{a_1\overline a_1}=t_1+\overline t_2, \ \ \ \
i\frac{a_2}{a_1\overline a_1}=t_2,
\endaligned
\end{equation}
which give, surprisingly, the parameters $t_1$ and $t_2$ as some weight zero expressions:
\begin{equation}
\label{t1,2-expression}
\boxed{t_1=\frac{a_3}{a_1^2\overline a_1}+2i\frac{\overline a_2}{a_1\overline a_1}, \ \ \ \ \ \ t_2=i\frac{a_2}{a_1\overline a_1}}.
\end{equation}
Putting these expressions in $\sf S$ and multiplying again the appearing fractional equations by some sufficient powers of $a_1$ and $\overline a_1$, then one finds $\sf S$ as a weighted homogeneous polynomial system {\it with no any parameter}. Except $a_2$ and $a_3$ that we already spent their associated equations \thetag{\ref{t1-2}} to find the expressions of the parameters $t_1$ and $t_2$, for each other involving group parameters $a_\bullet$ there exists one equation in $\sf S$ that expresses it in terms of some lower weight group parameters. Our next goal is to provide two more polynomial equations including $a_2$ and $a_3$ to recover this constraint.

 \subsubsection{Second part: structure equations of the weight $-\rho$}

One might be somehow surprised that so far we did not talk about the weight $-\rho$ structure equations. In fact, our trick was to retain them for our current aim of providing at least two more weighted homogeneous equations\footnote{Actually in CR dimension $1$, the reason of satisfying Beloshapka's maximum conjecture in the lengths $\rho\geqslant 3$ may refer to this part of our constructions. In fact, to provide two more weighted homogeneous equations for $a_2$ and $a_3$, we need some more structure equations than those of $d\Gamma_{2,3}, d\Gamma_{1,1}$ and $d\Gamma_{1,2}=\overline{d\Gamma_{1,1}}$. This means that we should at least have the next structure equation $d\Gamma_{3,4}$ which appears in the case of CR models which are of length $\rho\geqslant 3$.}. Notice that the method suggested above, can not be applied on a weight $-\rho$ structure equation:
\begin{equation}
\label{t1,2-a2,3}
d\Gamma_{\rho,i}=(p_i\alpha+q_i\overline\alpha)\wedge\Gamma_{\rho,i}+a_1^{p_i}\overline a_1^{q_i}\,d\sigma_{\rho,i}
\end{equation}
since it essentially does not contain any term $a_\bullet\,d\sigma_\bullet$ with $d\sigma_\bullet$ of the weight $-(\rho+1)$.
However, here we can think about picking up coefficients of the wedge products $\Gamma_{\rho,i}\wedge\Gamma_{1,t}$ from $d\Gamma_{\rho,i}$ for $t=1,2$. For this purpose, one notices that in the Darboux-Cartan structure of $d\sigma_{\rho,i}$, visible in the structure equation $d\Gamma_{\rho,i}$, only wedge products of the form $\sigma_{\rho-1,j}\wedge\sigma_{1,t}$ can make nonzero coefficients of $\Gamma_{\rho,i}\wedge\Gamma_{1,t}$. In order to find these coefficients and according to Lemma \ref{l-l-1}, if the coefficient of $d\sigma_{\rho,i}$ in the structure equation $d\Gamma_{\rho-1,j}$ is a (possibly zero) weight $\rho$ group parameter $a_{j_r}$, then the coefficient of $\Gamma_{\rho,i}$ in  $\sigma_{\rho-1,j}$ is some fraction of the form $-\frac{a_{j_r}}{a_1^\bullet\overline a_1^\bullet}$\,\,---\,\,notice that by considering the term $a_{j_r}d\sigma_{\rho,i}$ in the weight $-(\rho-1)$ structure equation $d\Gamma_{\rho-1,j}$, we find a weighted homogeneous equation of $\sf S$, expressing $a_{j_r}$ in terms of some lower weight group parameters. Then the coefficient of the sought wedge product $\Gamma_{\rho,i}\wedge\Gamma_{1,t}$ in $\sigma_{\rho-1,j}\wedge\sigma_{1,t}$ is the multiplication between the coefficient $-\frac{a_{j_r}}{a_1^\bullet\overline a_1^\bullet}$ of $\Gamma_{\rho,i}$ in $\sigma_{\rho-1,j}$ and the coefficient of $\Gamma_{1,t}$ in $\sigma_{1,t}$, which is $\frac{1}{a_1}$ where $t=1$ and $\frac{1}{\overline a_1}$ where $t=2$.  This implies that: (i)
 after absorption \thetag{\ref{replacement-maurer-cartan}} and equating to zero the coefficients of $\Gamma_{\rho,i}\wedge\Gamma_{1,1}$ and $\Gamma_{\rho,i}\wedge\Gamma_{1,2}$ in the structure equation $d\Gamma_{\rho,i}$, one finds two equations:
 \begin{equation}
 \label{a2-3-eq}
 \aligned
 \sum_{j_r}\,\frac{a_{j_r}}{a_1^\bullet\overline a_1^\bullet}+p_i\,t_1+q_i\,\overline t_2=0 \ \ \ \ {\rm and} \ \ \ \
  \sum_{j'_r}\,\frac{a_{j'_r}}{a_1^\bullet\overline a_1^\bullet}+q_i\,\overline t_1+p_i\,t_2=0,
 \endaligned
 \end{equation}
where according to \thetag{\ref{t1,2-expression}} they are actually two equations in terms of $a_2$, $a_3$ and some other weight $\rho$ group parameters $a_{j_r}$.
 (ii) In the system $\sf S$, one finds some polynomial equations which express $a_{j_r}$s and $a'_{j_r}$s in terms of some lower weight group parameters.

Surprisingly, Proposition \ref{prop-S} and equations \thetag{\ref{t1,2-expression}} imply that one can regard eventually the above two equations \thetag{\ref{a2-3-eq}} in terms of only $a_3, a_2, a_1$ and their conjugations. Now to finalize constructing the desired subsystem, it remains only to add these already found equations to $\sf S$.

Before attempt to solve the system $\sf S$, let us summarize our practical method of its construction. It is divided into the following two parts which should be performed after the absorption step \thetag{\ref{replacement-maurer-cartan}}:

\medskip
\noindent
{\bf Part I.} For each structure equation:
\[
\aligned
d\Gamma_{\ell,m}=(p_m\alpha+q_m\overline\alpha)\wedge\Gamma_{\ell,m}+\sum_{l\gneqq\ell}\,\delta_{i_t}\wedge\Gamma_{l,j}
+\sum_{l\geq\ell +2 }\,a_{j_n}d\sigma_{l,n}+\sum_{r}\,a_{j_r}d\sigma_{\ell+1,r}+a_1^{p_m}\overline a_1^{q_m}\,d\sigma_{\ell,m},
\endaligned
\]
with $\ell=1,\ldots,\rho-1$ and for each term $a_{j_{r_0}}d\sigma_{\ell+1,r_0}$ in its penultimate sum, equate to zero the coefficient of the wedge product $\Gamma_{\ell,i_j}\wedge\Gamma_{1,t_r}$, where $\sigma_{\ell,i_j}\wedge\sigma_{1,t_r}$ uniquely appears in the Darboux-Cartan structure of $d\sigma_{\ell+1, r_0}$ according to Lemma \ref{lem-unique-1-form}. The achieved equation belongs to $\sf S$.

\medskip
\noindent
{\bf Part II.} For each weight $-\rho$ structure equation:
\[
d\Gamma_{\rho,i}=(p_i\alpha+q_i\overline\alpha)\wedge\Gamma_{\rho,i}+a_1^{p_i}\overline a_1^{q_i}\,d\sigma_{\rho,i},
\]
equate to zero all coefficients of $\Gamma_{\rho,i}\wedge\Gamma_{1,t}$ for $t=1,2$ and add the achieved equations to $\sf S$.

\subsection{Solving the picked up subsystem}
\label{picking}

After constructing the weighted homogeneous polynomial system $\sf S$, now let us attempt to find the weighted projective variety ${\bf V}(\mathcal I)$ of the polynomial ideal $\mathcal I=\langle\sf S \rangle$\,\,---\,\,namely the solution of the system $\sf S$\,\,---\,\,in the weighted projective space $\mathbb P(1,2,3,\ldots)$ ({\it see e.g.} \cite{Alg-Geometry} for more details). Since the only weight $1$ group parameter $a_1$ is assumed to be nonzero, then this variety does not contain any point at the infinity surface $a_1=0$. Assume that $\mathcal I^{\sf aff}\subset\mathbb C[a_2,a_3,\ldots,a_r]$ is the affine ideal obtained as the dehomogenization of $\mathcal I$ by setting $a_1=1$. If $g$ is a weighted homogeneous polynomial in $\mathcal I$, then the following relation holds between it and its dehomogenization $g^{\sf deh}$ ({\it cf.} \cite[Theorem 5.16]{Alg-Geometry}):
\begin{equation}
\label{g-gdeh}
g(a_1,a_2,a_3,\ldots,a_r)=a_1^{\sf w-deg}\cdot g^{\sf deh}\big(\frac{a_2}{a_1^{[a_2]}}, \frac{a_3}{a_1^{[a_3]}},\ldots,\frac{a_r}{a_1^{[a_r]}}\big)
\end{equation}
where ${\sf w-deg}$ is the weight degree of the affine polynomial $g^{\sf deh}$. By Proposition \ref{prop-S} we can still state that
 associated with each group parameter $a_j$ visible in $\mathcal I^{\sf aff}$, there exists some (not necessarily weighted homogeneous, any more) polynomial in this ideal, expressed in terms of $a_j$ and some other group parameters (variables) of absolutely lower weights. Moreover, these polynomials are all linear (consider the equations of $\sf S$ after setting $a_1=1$ in \thetag{\ref{syst-ell}}, \thetag{\ref{t1,2-expression}}, \thetag{\ref{a2-3-eq}}).

This indicates that after selecting some appropriate order $\prec$ on the extant group parameters $a_\bullet$ enjoying the property that $a_i\prec a_j$ whenever $[a_i] < [a_j]$, then the affine ideal $\mathcal I^{\sf aff}$ is in fact in {\it Noether normal position} and according to the {\sl Finiteness Theorem} (\cite[Theorem 6 and Corollary 7, pp. 230-1]{Ideals-varieties}), the affine variety ${\bf V}(\mathcal I^{\sf aff})$ is zero dimensional containing just the origin $(0,0,\ldots,0)$.
Then according to the above equality \thetag{\ref{g-gdeh}}, one concludes that the weighted projective variety ${\bf V}(\mathcal I)$, or equivalently the solution set of the weighted homogeneous system $\sf S$, comprises some points of the concrete form:
\[
{\bf V}(\mathcal I)=\{(a_1,0,0,\ldots,0), \ \ \ \ \ a_1\neq 0\}.
\]
In other words, in the solution set of our weighted homogeneous system $\sf S$, all the group parameters visible in it\,\,---\,\,but not necessarily all the group parameters appearing in our ambiguity matrix\,\,---\,\,take the value zero, identically. In particular, the two fundamental group parameters $a_2$ and $a_3$ shall be zero. But, thanks to Lemma \ref{lem-a2-a3}, vanishing of these two group parameters is sufficient to assert that {\it all} the group parameters $a_j, j\neq 1$, appearing in the ambiguity matrix $\bf g$ should be normalized to zero;

\begin{Proposition}
\label{prop-aj=0}
After sufficient steps of applying absorption and normalization on the structure equations of the equivalence problem to a totally nondegenerate CR model $M_k$ of CR dimension $1$ and codimension $k$, all the appearing group parameters $a_j$ with $j=2,3,4,\ldots$ vanish, identically.
\end{Proposition}

This immediately results in the reduction of our ambiguity matrix group $G$ ({\it cf.} \thetag{\ref{ambiguity-coframe-k}}) to the simple diagonal matrix Lie group $G^{\sf red}$ comprising matrices of the form:
\begin{equation}
\label{G-red}
\footnotesize\aligned
{\bf g}^{\sf red}:=
\left(
  \begin{array}{cccc}
    a_1^{p}\overline a_1^{q} &  0 & \cdots & 0 \\
    \vdots &  \ddots & 0 & 0 \\
    0 & \cdots & \overline a_1 & 0 \\
    0 &  \cdots & 0 & a_1 \\
  \end{array}
\right).
\endaligned
\end{equation}
 Concerning the Maurer-Cartan matrix $\omega_{\tt MC}$ visible in \thetag{\ref{differentiation-2}}, all the Maurer-Cartan forms $\delta$ vanish identically and it reduces to a diagonal matrix with some combinations of the 1-forms $\alpha=\frac{da_1}{a_1}$ and its conjugation at its diagonal.  Finally, after  vanishing of the group parameters $a_2, a_3, \ldots$, then all torsion coefficients $T^i_{j,m}$ vanish identically except those which were constant from the beginning of construction;

\begin{Proposition}
\label{prop-struc-after-vanish}
After vanishing the group parameters $a_2, a_3, a_4,\ldots$, our structure equations convert into the simple {\sl constant type}:
\begin{equation}
\label{structure-equations-old-before-prolongation}
\aligned
d\,\Gamma_{\ell,i}&:=(p_i\,\alpha+q_i\,\overline\alpha)\wedge\Gamma_{\ell,i}
+\sum_{l+m=\ell\atop{j,n}}\,{\sf c}^{i}_{j,n}\,\Gamma_{l,j}\wedge\Gamma_{m,n} \ \ \ \ {\scriptstyle (\ell\,=\,1\,,\,\ldots\,,\,\rho\,,\, \ \ i\,=\,1\,,\,\ldots\,,\,2+k)}
\endaligned
\end{equation}
for some constant complex integers ${\sf c}^i_{j,n}$.
\end{Proposition}

\proof
According to \thetag{\ref{differentiation-2}}, our structure equations were originally of the form:
\[
d\Gamma_{\ell,i}=(p_i\alpha+q_i\overline\alpha)\wedge\Gamma_{\ell,i}+\zero{\sum_{l\gneqq\ell}\,\delta_{i_j}\wedge\Gamma_{l,j}}+
\zero{\sum_{l\gneqq\ell}\,a_{i_j}d\sigma_{l,j}}+a_1^{p_i}\overline a_1^{q_i}\,d\sigma_{\ell,i}.
\]
As mentioned, after vanishing of the group parameters $a_2, a_3, \ldots$  all the Maurer-Cartan forms $\delta_\bullet$ vanish identically and this kills the first sum $\sum_{l\gneqq\ell}\,\delta_{i_j}\wedge\Gamma_{l,j}$. For the second sum $\sum_{l\gneqq\ell}\,a_{i_j}d\sigma_{l,j}$ and according to Lemma \ref{lemma-weight-adsigma}, since all differentiations $d\sigma_{l,j}$ are of the weights $\lvertneqq -1$ (notice that here $l\gneqq\ell$ and $\ell\geqslant 1$) then all the group parameters $a_{i_j}$ are of the weights $\gneqq 1$ and hence none of them is $a_1$. This yields vanishing of this term, as well. Then, it suffices to consider the last term $a_1^{p_i}\overline a_1^{q_i}\,d\sigma_{\ell,i}$ of the above structure equation. According to the computed Darboux-Cartan structure in Proposition \ref{Darboux-Cartan} we have:
\[
d\sigma_{\ell,i}:=\sum_{\beta+\gamma=\ell}\,{\sf c}_{r,s}\,\sigma_{\beta,r}\wedge\sigma_{\gamma,s}.
\]
On the other hand, our inverse matrix ${\bf g}^{-1}$ is now converted to the simple form:
\[
\footnotesize
({\bf g^{\sf red}})^{-1}
=
\left(
  \begin{array}{cccc}
    \frac{1}{a_1^{p}\overline a_1^{q}} &  0 & \cdots & 0 \\
    \vdots &  \ddots & 0 & 0 \\
    0 &  \cdots & \frac{1}{\overline a_1} & 0 \\
    0 &  \cdots & 0 & \frac{1}{a_1} \\
  \end{array}
\right)
\]
which through the equality $\Sigma=({\bf g}^{\sf red})^{-1}\cdot\Gamma$, it results that:
\[
\sigma_{\beta,r}\wedge\sigma_{\gamma,s}=\frac{1}{a_1^{m_r}\overline a_1^{n_s}}\,\Gamma_{\beta,r}\wedge\Gamma_{\gamma,s},
\]
for some constant integers $m_r$ and $n_s$. Then, the last term $a_1^{p_i}\overline a_1^{q_i}\,d\sigma_{\ell,i}$ can be brought into a combination as:
\[
a_1^{p_i}\overline a_1^{q_i}\,d\sigma_{\ell,i}:=\sum_{\beta+\gamma=\ell}\,{\sf c}_{r,s}\frac{a_1^{p_i}\overline a_1^{q_i}}{a_1^{m_r}\overline a_1^{n_s}}\,\Gamma_{\beta,r}\wedge\Gamma_{\gamma,s}.
\]
On the other hand, these coefficients ${\sf c}_{r,s}\frac{a_1^{p_i}\overline a_1^{q_i}}{a_1^{m_r}\overline a_1^{n_s}}$ are in fact the only remained torsion coefficients $T^i_{rs}$ of the wedge products $\Gamma_{\beta,r}\wedge\Gamma_{\gamma,s}$, in the structure equation $d\Gamma_{\ell,i}$ and hence according to Proposition \ref{prop-torsion-weight-zero}, are of the weight zero. Since they involve just weight one group parameters $a_1$ and $\overline a_1$ then, after simplifications if necessary, they will be either some constants or some fractions  of the form:
\[
{\sf c}_{r,s}\frac{a_1^i}{\overline a_1^i} \ \ \ \ \textrm{or} \ \ \ \ {\sf c}_{r,s}\frac{\overline a_1^i}{a_1^i}.
\]
Consequently, our structure equation $d\Gamma_{\ell,i}$ is now converted into the form:
\[
\aligned
d\Gamma_{\ell,i}&=(p_i\alpha+q_i\overline\alpha)\wedge\Gamma_{\ell,i}+\sum_{\beta'+\gamma'=\ell}\,{\sf c}_{r',s'}\,\Gamma_{\beta',r'}\wedge\Gamma_{\gamma',s'}+
\\
&\ \ \ \ \ \ \ \ \ \ \ \ \ +\sum_{\beta+\gamma=\ell}\,{\sf c}_{r,s}\frac{a_1^i}{\overline a_1^i}\,\Gamma_{\beta,r}\wedge\Gamma_{\gamma,s}+\sum_{\beta+\gamma=\ell}\,{\sf c}_{r,s}\frac{\overline a_1^j}{a_1^j}\,\Gamma_{\beta,r}\wedge\Gamma_{\gamma,s}.
\endaligned
\]
All the appearing $\beta$s and $\gamma$s in this equation are absolutely less than $\ell$, whence in the case that one ${\sf c}_{r,s}$ is nonzero then the torsion coefficient $T^i_{r,s}={\sf c}_{r,s}\frac{a_1^i}{\overline a_1^i}$ or $T^i_{r,s}={\sf c}_{r,s}\frac{\overline a_1^i}{ a_1^i}$ of $\Gamma_{\beta,r}\wedge\Gamma_{\gamma,s}$ is essential and can be plainly normalized to some constant, say ${\sf c}_{r,s}$, by normalizing  $\frac{a_1}{\overline a_1}=1$, {\it i.e.} by considering the only remained parameter $a_1$ as {\it real}. Then all powers of $\frac{a_1}{\overline a_1}$ will be equal to $1$ and consequently, we receive finally just some constant coefficients of these remaining wedge products.
\endproof

What mentioned at the end of the above proof also demonstrates the normalization of the only remained group parameter $a_1$. Accordingly, this parameter is never normalizable in the case that after vanishing the group parameters $a_2, a_3, \ldots$, all the torsion coefficients of the structure equations are constant. Otherwise, $a_1$ will be normalized just to a real group parameter.

\begin{Corollary}
\label{cor-a1}
There are two possibility for the normalization of the only remained group parameter $a_1$. It is either normalizable to a real group parameter or it is never normalizable. The reduced structure group $G^{\sf red}$ (cf. \thetag{\ref{G-red}}) is of real dimension $1$ in the former case and $2$ in the latter.
\end{Corollary}

For instance, one observes in \cite{5-cubic} that in the case of $M_3$, the group parameter $a_1$ is never normalizable while, in contrary, for $M_4$ it is normalizable to a real group parameter as is shown in \cite{SCM-Moduli}.

\subsection{Prolongation}

To continue toward Cartan's approach of solving equivalence problems and after applying sufficient absorption-normalization steps, now one has to start the {\it prolongation} step. The main result behind this step is \cite[Proposition 12.1]{Olver-1995}. This Proposition permits us to substitute the current equivalence problem to the $(2+k)$-dimensional CR model $M_k$ by that of the $(3+k)$ or $(4+k)$-dimensional prolonged space $M^{\sf pr}:=M_k\times G^{\sf red}$. For this, we have to add the remaining Maurer-Cartan forms $\alpha$ and $\overline\alpha$ to the original lifted coframe $\Gamma$ and consider $(\Gamma_{\rho,2+k},\ldots,\Gamma_{1,1},\alpha,\overline\alpha)$ as the new lifted coframe associated with this prolonged space. In the case that $a_1$ is normalizable to a real group parameter, then of course we have $\alpha=\overline\alpha$. Constructing the associated structure equations to this new problem is easy, just adding $d\alpha=d\big(\frac{d\,a_1}{a_1}\big)=0$ to the former ones. Then, the final structure equations of our new equivalence problem to the prolonged space $M^{\sf pr}$ take the following $\{e\}$-{\it structure} constant type:
 \begin{equation}
\label{structure-equations-after-prolongation}
\aligned
\left [
\begin{array}{l}
d\,\Gamma_{\ell,i}=(p_i\,\alpha+q_i\,\overline\alpha)\wedge\Gamma_{\ell,i}
+\sum_{\ell_1+\ell_2=\ell}\,{\sf c}^{i}_{j,n}\,\Gamma_{\ell_1,j}\wedge\Gamma_{\ell_2,n} \ \ \ \ {\scriptstyle (\ell\,=\,1\,,\,\ldots\,,\,\rho, \ \ i\,=\,1\,,\,\ldots\,,\,2+k)},
\\
d\alpha=0,
\\
d\overline\alpha=0.
\end{array}
\right.
\endaligned
\end{equation}

Then we have arrived at the stage of stating the main result of this paper;

\begin{Theorem}
\label{theorem-main}
The biholomorphic equivalence problem to a $(2+k)$-dimensional real analytic totally nondegenerate CR model $M_k\subset\mathbb C^{1+k}$ of codimension $k$ is reducible to some absolute parallelisms, namely to some certain $\{e\}$-structures on prolonged manifolds $M_k\times G^{\sf red}$ of real dimensions either $3+k$ or $4+k$.
\end{Theorem}

\medskip\noindent
{\bf Weight assignment.} We assign naturally\footnote{Notice that the exterior differentiation $d\alpha$ took the value zero, exactly as constant functions.} the weight zero to the new lifted 1-forms $\alpha$ and $\overline\alpha$.

\section{Proof of Beloshapka's maximum conjecture}
\label{section-proof}

After providing key results in the previous section, now we are ready to prove Beloshapka's maximum conjecture \ref{conjecture} in CR dimension one. As we saw, the equivalence problem to a certain CR model $M_k$ converted finally to that of the prolonged space $M^{\sf pr}$ with the final constant type structure equations \thetag{\ref{structure-equations-after-prolongation}}.
According to \cite[Theorem 8.16]{Olver-1995}, if the final structure equations of an equivalence problem to an $r$-dimensional smooth manifold $M$ equipped with some lifted coframe $\{\gamma^1,\ldots,\gamma^r\}$ is of the constant type:

\[
d\gamma^k
=
\sum_{1\leqslant i<j\leqslant r}\,
c^k_{ij}\,\gamma^i\wedge\gamma^j
\ \ \ \ \ \ \ \ \ \ \ \ \
{\scriptstyle{(k\,=\,1\,\cdots\,r),}}
\]
then $M$ is (locally) diffeomorphic to an $r$-dimensional Lie group $\sf G$ corresponding to the Lie algebra $\frak g$ with the basis elements $\{{\sf v}_1,\ldots, {\sf v}_r\}$ and enjoying the so-called {\it structure constants}:
\[
\big[{\sf v}_i,{\sf v}_j\big]
=
-\sum_{k=1}^r\,c^k_{ij}\,{\sf v}_k
\ \ \ \ \ \ \ \ \ \ \ \ \
{\scriptstyle{(1\,\leqslant\,i\,<\,j\,\leqslant\,r)}}.
\]

Accordingly, let us try to find the Lie algebra $\frak g$ corresponding to the constant structure equations \thetag{\ref{structure-equations-after-prolongation}}. At first, we associate to each lifted 1-form $\Gamma_{\ell,i}$ of $M^{\sf pr}$ the basis element ${\sf v}_{\ell,i}$ of $\frak g$. For the new appearing lifted 1-forms $\alpha$ and $\overline\alpha$, let us associate ${\sf v}_0$ and ${\sf v}_{\overline 0}$. Of course, if the real part of $a_1$ is normalizable, then we dispense with ${\sf v}_{\overline 0}$ since in this case we have $\alpha=\overline\alpha$ and hence our desired Lie algebra $\frak g$ is of dimension either $3+k$ or $4+k$, depending on the normalization of $a_1$. Assign naturally the weight $-\ell$ to each basis element ${\sf v}_{\ell,i}$ and the weight zero to ${\sf v}_0$ and ${\sf v}_{\overline 0}$. In particular, because we do not see any wedge product $\alpha\wedge\overline\alpha$ among the structure equations \thetag{\ref{structure-equations-after-prolongation}} then, $[{\sf v}_0,{\sf v}_{\overline 0}]=0$. This indicates that $\{{\sf v}_0,{\sf v}_{\overline 0}\}$ generates an {\it Abelian} subalgebra of $\frak g$.

Each structure equation $d\Gamma_{\ell,i}$ in \thetag{\ref{structure-equations-after-prolongation}} is some constant combination of the wedge products between lifted 1-forms for which the sum of their weights is exactly $-\ell$. Thus, the Lie bracket between two weight  $-\ell_1$ and $-\ell_2$ basis elements ${\sf v}_\bullet$ of $\frak g$ will be some constant combination of its weight $-(\ell_1+\ell_2)$ basis elements. Consequently, we have the following interesting result;

\begin{Proposition}
\label{prop-graded-algebra}
Let $\frak g_{-\ell}$ be the $\mathbb C$-vector space generated by all basis elements ${\sf v}_{\ell,i}$ of the weight $-\ell$ and let $\frak g_0$ is the  Abelian subalgebra of $\frak g$ generated by ${\sf v}_0$ and ${\sf v}_{\overline 0}$. Then, the Lie algebra $\frak g$ associated with the final structure equations \thetag{\ref{structure-equations-after-prolongation}} is graded of the form:
\[
\frak g:=\underbrace{\frak g_{-\rho}\oplus\frak g_{-(\rho-1)}\oplus\ldots\oplus\frak g_{-1}}_{\frak g_-}\oplus\frak g_{0}
\]
satisfying $[\frak g_{-\ell_1},\frak g_{-\ell_2}]=\frak g_{-(\ell_1+\ell_2)}$. In this case, $\frak g_-$ is $(2+k)$-dimensional and $\frak g_0$ is of dimension either $1$ or $2$.
\end{Proposition}

On the other hand, Corollary 14.20 of \cite{Olver-1995} says that this Lie algebra $\frak g$ is in fact the symmetry Lie algebra of the prolonged space $M^{\sf pr}=M_k\times G^{\sf red}$ with respect to its coframe $(\Gamma_{1,1}, \ldots, \Gamma_{\rho,2+k}, \alpha,\overline\alpha)$; that is the Lie algebra associated with the Lie group $\sf G$ of self-equivalences $\Phi:M^{\sf pr}\rightarrow M^{\sf pr}$, satisfying $\Phi^\ast(\theta)=\theta$ for $\theta=\Gamma_{1,1}, \ldots, \Gamma_{\rho,2+k}, \alpha,\overline\alpha$. But, according to \cite[Proposition 12.1]{Olver-1995} and its proof, $\sf G$ can be identified with the CR symmetry Lie group ${\sf Aut}_{CR}(M)$ of biholomorphic maps $h:M_k\rightarrow M_k$, hence:
\[
\frak{aut}_{CR}(M_k)\cong\frak g.
\]
Consequently similar to $\frak g$, the Lie algebra $\frak{aut}_{CR}(M_k)$ will be graded without any positive component in its gradation, as was conjecture by Beloshapka.

\begin{Theorem}
{\sl (Beloshapka's maximum conjecture in CR dimension one)}. The Lie algebra $\frak{aut}_{CR}(M_k)$ associated with a Beloshapka's real analytic totally nondegenerate  CR model $M_k$ of CR dimension $1$, codimension $k$ and length $\rho\geq 3$\,\,---\,\,or equivalently of codimension $k\geq 2$\,\,---\,\,contains no any homogeneous component of absolutely positive homogeneity. In other words, such CR model has rigidity. Moreover, this Lie algebra is of dimension either $3+k$ or $4+k$.
\qed
\end{Theorem}

\appendix
\section{An example in the length four}
\label{appendix}

By way of illustration the method introduced in Section \ref{section-picking-system}, in this appendix we consider the biholomorphic equivalence problem to the $8$-dimensional, length $\rho=4$ CR model $M_6\subset\mathbb C^7$ represented as the graph of six defining polynomials:
\[
\aligned
w_1-\overline w_1&=2i\,z\overline z,
\\
w_2-\overline w_2&=2i\big(z^2\overline z+z\overline z^2\big), \ \ \ \ w_3-\overline w_3=2\big(z^2\overline z-z\overline z^2\big),
\\
w_4-\overline w_4&=2i\big(z^3\overline z+z\overline z^3\big), \ \ \  w_5-\overline w_5=2\big(z^3\overline z-z\overline z^3\big), \ \ \ w_6-\overline w_6=2i\,z^2\overline z^2.
\endaligned
\]
The assigned weights to the extant complex variables are:
\[
[z]=1, \ \ \ \ [w_1]=2, \ \ \ \ [w_2]=[w_3]=3, \ \ \ \ [w_4]=[w_5]=[w_6]=4.
\]

Saving the space, we do not present the intermediate calculations. According to our computations, our initial frame contains eight vector fields of various lengths $-1,\ldots,-4$:
\[
\aligned
\mathcal L&:=\mathcal L_{1,1}, \ \ \  \ \ \ \overline{\mathcal L}:=\mathcal L_{1,2},
\\
 \mathcal T&:=\mathcal L_{2,3}=i[\mathcal L, \overline{\mathcal L}],
 \\
 \mathcal S&:=\mathcal L_{3,4}=[\mathcal L, \mathcal T], \ \ \ \ \ \ \overline{\mathcal S}:=\mathcal L_{3,5}=[\overline{\mathcal L}, \mathcal T],
 \\
 \mathcal U&:=\mathcal L_{4,6}=[\mathcal L, \mathcal S], \ \ \ \ \ \ \overline{\mathcal U}:=\mathcal L_{4,7}=[\overline{\mathcal L}, \overline{\mathcal S}], \ \ \ \ \ \ \mathcal V:=\mathcal L_{4,8}=[\mathcal L, \overline{\mathcal S}]=[\overline{\mathcal L}, \mathcal S].
\endaligned
\]
The other Lie brackets between these eight initial vector fields are all zero. Assume that:
\[
\Sigma:=\big(\underbrace{\nu_0,\mu_0,\overline\mu_0}_\textrm{ weight -4}, \underbrace{\sigma_0, \overline\sigma_0}_\textrm{ weight -3}, \underbrace{\rho_0}_\textrm{ weight -2}, \underbrace{\zeta_0, \overline\zeta_0}_\textrm{ weight -1}\big)^t \,\,\,\, \textrm{ is the dual coframe of}  \,\,\,\, \big(\mathcal V, \mathcal U, \overline{\mathcal U}, \mathcal S, \overline{\mathcal S}, \mathcal T, \mathcal L, \overline{\mathcal L}\big)^t.
\]
Then the associated Darboux-Cartan structure to this coframe is:
\[
\aligned
d\nu_0&=\overline\sigma_0\wedge\zeta_0+\sigma_0\wedge\overline\zeta_0, \ \ \ \ d\mu_0=\sigma_0\wedge\zeta_0, \ \ \ \ d\overline\mu_0=\overline\sigma_0\wedge\overline\zeta_0
\\
d\sigma_0&=\rho_0\wedge\zeta_0, \ \ \ \ \ \ \ \ \ \ \ \ \ \ \ \ \  \ \ \ \ d\overline\sigma_0=\rho_0\wedge\overline\zeta_0,
\\
d\rho_0&=i\zeta_0\wedge\overline\zeta_0, \ \ \ \ \ \ \ \ \ \ \ \ \ \ \ \ \ \  \ \ d\zeta_0=0, \ \ \ \ \ \ \ \ \ \ \ \ \ \ \   d\overline\zeta_0=0.
\endaligned
\]
Assuming $\Gamma:=(\nu,\mu,\overline\mu, \sigma, \overline\sigma, \rho, \zeta, \overline\zeta\big)^t$ as the associated lifted coframe, then our computation brings the ambiguity $8\times 8$ invertible matrix of the biholomorphic equivalence problem to $M_6$ as:
\begin{equation}
\label{a6}
\footnotesize\aligned
\Gamma
=
\underbrace{\left(
  \begin{array}{cccccccc}
    a_1^2\overline a_1^2 & 0 & 0 & 0 & 0 & 0 & 0 & 0 \\
    0 & a_1^3\overline a_1 & 0 & 0 & 0 & 0 & 0 & 0 \\
    0 & 0 & a_1\overline a_1^3 & 0 & 0 & 0 & 0 & 0 \\
    a_{13} & a_6 & 0 & a_1^2\overline a_1 & 0 & 0 & 0 & 0 \\
    \overline a_{13} & 0 & \overline a_6 & 0 & a_1\overline a_1^2 & 0 & 0 & 0 \\
    a_{11} & a_7 & \overline a_7 & a_3 & \overline a_3 & a_1\overline a_1 & 0 & 0 \\
    a_{12} & a_8 & \overline a_9 & a_4 & a_5 & a_2 & a_1 & 0 \\
    \overline a_{12} & a_9 & \overline a_8 & \overline a_5 & \overline a_4 & \overline a_2 & 0 & \overline a_1 \\
  \end{array}
\right)}_{\bf g}
\cdot
\Sigma,
\endaligned
\end{equation}
with the assigned weights:
\[
[a_1]=1, \ \ \ [a_2]=2, \ \ \ [a_3]=[a_4]=[a_5]=3, \ \ \ [a_6]=\ldots=[a_{13}]=4.
\]

By computing the (somehow big) inverse matrix ${\bf g}^{-1}$, one can check also the assertion of some results like Lemmas \ref{lemma-weight-g-1} and \ref{sigma-Gamma-expression}. Also, our Maurer-Cartan matrix is of the form:

\[
\footnotesize\aligned
\omega_{\tt MC}:=\left(
  \begin{array}{cccccccc}
    2\alpha+2\overline \alpha & 0 & 0 & 0 & 0 & 0 & 0 & 0 \\
    0 & 3\alpha+\overline\alpha & 0 & 0 & 0 & 0 & 0 & 0 \\
    0 & 0 & \alpha+3\overline \alpha & 0 & 0 & 0 & 0 & 0 \\
    \delta_{13} & \delta_6 & 0 & 2\alpha+\overline \alpha & 0 & 0 & 0 & 0 \\
    \overline \delta_{13} & 0 & \overline \delta_6 & 0 & \alpha+2\overline \alpha& 0 & 0 & 0 \\
    \delta_{11} & \delta_7 & \overline \delta_7 & \delta_3 & \overline \delta_3 & \alpha+\overline \alpha & 0 & 0 \\
    \delta_{12} & \delta_8 & \overline \delta_9 & \delta_4 & \delta_5 & \delta_2 & \alpha& 0 \\
    \overline \delta_{12} & \delta_9 & \overline \delta_8 & \overline \delta_5 & \overline \delta_4 & \overline \delta_2 & 0 & \overline \alpha \\
  \end{array}
\right), \ \ \ \ \ {\rm with} \ \ \ \ \ \alpha=\frac{da_1}{a_1}.
\endaligned
\]
Then, our structure equations will be of the form\,\,---\,\,we abbreviate the superfluous combinations of the wedge products $\delta_j\wedge\bullet$ just by some $"\cdots "$ since they will not play any important role:
\begin{equation}
\label{struc-k=6}
\aligned
d\nu&=(2\alpha+2\overline\alpha)\wedge\nu+a_1^2\overline a_1^2\,d\nu_0,
\\
d\mu&=(3\alpha+\overline\alpha)\wedge\mu+a_1^3\overline a_1\,d\mu_0,
\\
d\sigma&=\cdots+(2\alpha+\overline\alpha)\wedge\sigma+a_{13}\,d\nu_0+a_6\,d\mu_0+a_1^2\overline a_1\,d\sigma_0,
\\
d\rho&=\cdots+(\alpha+\overline\alpha)\wedge\rho+a_{11}\,d\nu_0+a_7\,d\mu_0+\overline a_7\,d\overline\mu_0+a_3\,d\sigma_0+\overline a_3\,d\overline\sigma_0+a_1\overline a_1\,d\rho_0,
\\
d\zeta&=\cdots+\alpha\wedge\zeta+a_{12}\,d\nu_0+a_8\,d\mu_0+\overline a_9\,d\overline\mu_0+a_4\,d\sigma_0+a_5\,d\overline\sigma_0+a_2\,d\rho_0+a_1\,d\zeta_0.
\endaligned
\end{equation}

Now, let us proceed as subsection \ref{picking} to pick the appropriate weighted homogeneous system $\sf S$. To do it and as is the method of absorption-normalization step, first we apply the substitutions:
 \[
 \aligned
 \alpha&\mapsto\alpha+t_8\,\nu+t_7\,\mu+\ldots+t_2\,\overline\zeta+t_1\,\zeta,
 \\
  \delta_j&\mapsto\delta_j+s^j_8\,\nu+s^j_7\,\mu+\ldots+s^j_2\,\overline\zeta+s^j_1\,\zeta, \ \ \ \ \ \ \ {\scriptstyle (j\,=\,2\,,\,\ldots\,,\,13)}
 \endaligned
 \]
 on the above structure equations. According to our proposed method of constructing $\sf S$, in the minimum weight $-4$ structure equations $d\nu$ and $d\mu$, we have to compute the coefficients of $\nu\wedge\{\zeta,\overline\zeta\}$ and $\mu\wedge\{\zeta,\overline\zeta\}$, respectively. Moreover, in the weight $-3$ structure equation $d\sigma$, we should pick up the coefficients of $\sigma\wedge\{\zeta,\overline\zeta\}$ since $\sigma_0\wedge\zeta_0$ and $\sigma_0\wedge\overline\zeta_0$ uniquely appear in the Darboux-Cartan structure of the only extant length $-4$ differentiations $d\mu_0$ and $d\nu_0$ visible in this structure equation. Similarly, in the lengths $-2$ and $-1$ structure equations $d\rho$ and $d\zeta$, we should pick the coefficients of $\rho\wedge\{\zeta,\overline\zeta\}$ and $\zeta\wedge\overline\zeta$, respectively. Equating these coefficients to zero gives the following equations:
 \[
\footnotesize \aligned
 {\sf S}_{d\nu}&:=\big\{-\frac{\overline a_{13}}{a_1^2\overline a_1^2}=2t_1+2\overline t_2\big\}, \ \ \ \ \  \ \ \ \ \ \ \ \ \ \  \ \ \ \  \ \ \ \ \ \  {\sf S}_{d\mu}:=\big\{-\frac{a_6}{a_1^3\overline a_1}=3t_1+\overline t_2, \ \ \ 0=\overline t_1+3t_2\big\},
 \\
  {\sf S}_{d\sigma}&:=\big\{\frac{a_6}{a_1^3\overline a_1}-\frac{a_3}{a_1^2\overline a_1}=2t_1+\overline t_2, \ \ \ \frac{\overline a_{13}}{a_1^2\overline a_1^2}=\overline t_1+2t_2\big\},
  \ \ \ \ \ \
  {\sf S}_{d\rho}:=\big\{\frac{a_3}{a_1^2\overline a_1}+i\frac{\overline a_2}{a_1\overline a_1}=t_1+\overline t_2\big\},
  \ \ \ \ \ \   {\sf S}_{d\zeta}:=\big\{i\frac{ a_2}{a_1\overline a_1}=t_2\big\},
 \endaligned
 \]
where $\sf S$ is the union of them. Putting the obtained expressions of the parameters $t_1$ and $t_2$ into these equations and multiplying them by sufficient powers of $a_1$ and $\overline a_1$, one finds the following weighted homogeneous system:
 \[
 \aligned
 {\sf S}:=\bigg\{\overline a_{13}+2\,\overline a_1 a_3+2i\,a_1\overline a_1 \overline a_2=0, & \ \ \ \ a_6+3\,a_1a_3+5i\,a_1^2\overline a_2=0, \ \ \ \ \overline a_3+i\,\overline a_1a_2=0,
 \\
 & \ \ \ \  a_6-3\,a_1a_3-3i\,a_1^2\overline a_2=0, \ \ \ \ \overline a_{13}- a_1\,\overline a_3=0\bigg\}.
 \endaligned
 \]
 Either by hand or by means of some computer softwares, one versifies that the solution of this system is nothing but $a_2=a_3=a_6=a_{13}\equiv 0$, which immediately implies vanishing of all the group parameters $a_2,a_3,a_4,\ldots,a_{13}$. Our computations shows that here $a_1$ is not normalizable. Applying these results and after one prolongation, the first structure equations \thetag{\ref{struc-k=6}}  converts to the simple constant form:
\[
\aligned
d\nu&=(2\alpha+2\overline\alpha)\wedge\nu+\overline\sigma\wedge\zeta+\sigma\wedge\overline\zeta,
\\
d\mu&=(3\alpha+\overline\alpha)\wedge\mu+\sigma\wedge\zeta,
\\
d\sigma&=(2\alpha+\overline\alpha)\wedge\sigma+\rho\wedge\zeta,
\\
d\rho&=(\alpha+\overline\alpha)\wedge\rho+i\,\zeta\wedge\overline\zeta,
\\
d\zeta&=\alpha\wedge\zeta
\\
d\alpha&=0.
\endaligned
\]

\begin{Proposition}
The Lie algebra $\frak g$ associated with the above structure equations is $10$-dimensional with the basis $\{{\sf v}^{\nu}, {\sf v}^{\mu}, {\sf v}^{\overline\mu}, {\sf v}^{\sigma}, {\sf v}^{\overline\sigma}, {\sf v}^{\rho}, {\sf v}^{\zeta}, {\sf v}^{\overline\zeta}, {\sf v}^{\alpha}, {\sf v}^{\overline\alpha}\}$ and with the Lie brackets, displayed in the following table:
\medskip
\begin{center}
\footnotesize
\begin{tabular} [t] { c | c c c c c c c c c c}
& ${\sf v}^{\nu}$ & ${\sf v}^{\mu}$ & ${\sf v}^{\overline\mu}$ & ${\sf v}^{\sigma}$ & ${\sf v}^{\overline\sigma}$ & ${\sf v}^{\rho}$ & ${\sf v}^{\zeta}$ & ${\sf v}^{\overline\zeta}$ & ${\sf v}^{\alpha}$ & ${\sf v}^{\overline\alpha}$
\\
\hline
${\sf v}^{\nu}$ & $0$ & $0$ & $0$ & $0$ & $0$ & $0$ & $0$ & $0$ & $2{\sf v}^{\nu}$ & $2{\sf v}^{\nu}$
\\
${\sf v}^{\mu}$ & $*$ & $0$ & $0$ & $0$ & $0$ & $0$ & $0$ & $0$ & $3{\sf v}^{\mu}$ & ${\sf v}^{\mu}$
\\
${\sf v}^{\overline\mu}$ & $*$ & $*$ & $0$ & $0$ & $0$ & $0$ & $0$ & $0$ & ${\sf v}^{\overline\mu}$ & $3{\sf v}^{\overline\mu}$
\\
${\sf v}^{\sigma}$ & $*$ & $*$ & $*$ & $0$ & $0$ & $0$ & $-{\sf v}^\mu$ & $-{\sf v}^\nu$ & $2{\sf v}^{\sigma}$ & ${\sf v}^{\sigma}$
\\
${\sf v}^{\overline\sigma}$ & $*$ & $*$ & $*$ & $*$ & $0$ & $0$ & $-{\sf v}^\nu$ & $-{\sf v}^{\overline\mu}$ & ${\sf v}^{\overline\sigma}$ & $2{\sf v}^{\overline\sigma}$
\\
${\sf v}^{\rho}$ & $*$ & $*$ & $*$ & $*$ & $*$ & $0$ & $-{\sf v}^\sigma$ & $-{\sf v}^{\overline\sigma}$ & ${\sf v}^{\rho}$ & ${\sf v}^{\rho}$
\\
${\sf v}^{\zeta}$ &  $*$ & $*$ & $*$ & $*$ & $*$ & $*$ & $0$ & $-i{\sf v}^{\rho}$ & ${\sf v}^{\zeta}$ & $0$
\\
${\sf v}^{\overline\zeta}$ &  $*$ & $*$ & $*$ & $*$ & $*$ & $*$ & $*$ & $0$ & $0$ & ${\sf v}^{\overline\zeta}$
\\
${\sf v}^{\alpha}$ &  $*$ & $*$ & $*$ & $*$ & $*$ & $*$ & $*$ & $*$ & $0$ & $0$
\\
${\sf v}^{\overline\alpha}$ &  $*$ & $*$ & $*$ & $*$ & $*$ & $*$ & $*$ & $*$ & $*$ & $0$
\end{tabular}
\end{center}
This Lie algebra, which is isomorphic to $\frak{aut}_{CR}(M_6)$, is graded of the form:
\[
\frak g:=\frak g_{-4}\oplus\frak g_{-3}\oplus\frak g_{-2}\oplus\frak g_{-1}\oplus\frak g_0,
\]
with $\frak g_{-4}=\langle{\sf v}^{\nu}, {\sf v}^{\mu}, {\sf v}^{\overline\mu}\rangle$, with $\frak g_{-3}=\langle {\sf v}^{\sigma}, {\sf v}^{\overline\sigma}\rangle$, with $\frak g_{-2}=\langle {\sf v}^{\rho}\rangle$, with $\frak g_{-1}=\langle {\sf v}^{\zeta}, {\sf v}^{\overline\zeta}\rangle$ and with $\frak g_0=\langle {\sf v}^{\alpha}, {\sf v}^{\alpha}\rangle$.
\end{Proposition}

\subsection*{Acknowledgment}

The author expresses his sincere thanks to Jo\"el Merker and Amir Hashemi for their
helpful comments and discussions during the preparation of this paper. The research of the author was supported in part by a grant from IPM, no. 95510427.

\end{document}